\documentclass[11pt]{amsart}

\usepackage[utf8]{inputenc}
\usepackage[T1]{fontenc}
\usepackage{lmodern}

\usepackage{amssymb}
\usepackage{mathtools}
\usepackage{amscd}
\usepackage[mathscr]{eucal}
\usepackage{enumitem}
\usepackage[hmargin=1in,vmargin=1in]{geometry}

\usepackage{graphicx}
\usepackage{xcolor}
\usepackage[roman]{sublabel}

\usepackage[labelfont=bf,font=small]{caption}
\usepackage{subcaption}

\usepackage[noadjust]{cite}
\usepackage[foot]{amsaddr}

\usepackage{placeins}
\usepackage{bm}
\usepackage{siunitx}

\definecolor{cite_purple}{RGB}{128,9,158}
\definecolor{cite_blue}{RGB}{2,95,176}
\definecolor{link_red}{rgb}{0.7,0,0}
\usepackage[colorlinks=true,
            citecolor=cite_purple,
            linkcolor=link_red,
            urlcolor=cite_blue]{hyperref}

\newlength{\arrayrulewidthOriginal}
\newcommand{\Cline}[2]{%
  \noalign{\global\setlength{\arrayrulewidthOriginal}{\arrayrulewidth}}%
  \noalign{\global\setlength{\arrayrulewidth}{#1}}\cline{#2}%
  \noalign{\global\setlength{\arrayrulewidth}{\arrayrulewidthOriginal}}}

\usepackage{ottmacros-shear-chaos-rev2}

\DeclareMathOperator{\Var}{Var}

\theoremstyle{plain}

\theoremstyle{definition}

\begin{document}

\title{Delay-induced uncertainty for a paradigmatic glucose-insulin model}
% \title{Reliability failure in physiological systems}

\author{Bhargav Karamched$^{1}$}
\address{$^{1}$Department of Mathematics and Institute of Molecular Biophysics, Florida State University, Tallahassee, Florida, 32306, USA}

\author{George Hripcsak$^{2}$}
\address{$^{2}$Department of Biomedical Informatics, Columbia University, New York, New York, 10032, USA}

\author{David Albers$^{3, 2}$}
\address{$^{3}$Section of Informatics and Data Science, Department of Pediatrics, Department of Biomedical Engineering, Department of Biostatistics and Informatics, University of Colorado Anschutz Medical Campus, Aurora, Colorado, 80045, USA}

\author{William Ott$^{*,4}$}
\address{$^{4}$Department of Mathematics, University of Houston, Houston, Texas, 77204, USA}
\thanks{$^{*}$\href{mailto:ott@math.uh.edu}{ott@math.uh.edu} (Corresponding author)}

% \urladdr[William Ott]{http://www.math.uh.edu/$\sim$ott/}

\keywords{Delay-induced uncertainty, glucose-insulin dynamics, Lyapunov exponent, reliability failure, shear, strange attractor, supercritical Hopf bifurcation, theory of rank one maps}

\subjclass[2010]{37D25, 37D45, 37G35, 37N25, 92C30}

\date{\today}

\begin{abstract}
Medical practice in the intensive care unit is based on the assumption that physiological systems such as the human glucose-insulin system are predictable.
We demonstrate that delay within the glucose-insulin system can induce sustained temporal chaos, rendering the system unpredictable.
Specifically, we exhibit such chaos for the Ultradian glucose-insulin model.
This well-validated, finite-dimensional model represents feedback delay as a three-stage filter.
Using the theory of rank one maps from smooth dynamical systems, we precisely explain the nature of the resulting delay-induced uncertainty (DIU).
We develop a framework one may use to diagnose DIU in a general oscillatory dynamical system.
For infinite-dimensional delay systems, no analog of the theory of rank one maps exists.
Nevertheless, we show that the geometric principles encoded in our DIU framework apply to such systems by exhibiting sustained temporal chaos for a linear shear flow.
Our results are potentially broadly applicable because delay is ubiquitous throughout mathematical physiology.
\end{abstract}

\maketitle

\raggedbottom
\thispagestyle{empty}

{\bfseries
We introduce a novel route through which delay causes oscillatory dynamical systems to exhibit sustained temporal chaos.
We precisely explain the nature of the resulting delay-induced uncertainty (DIU).
We show that DIU occurs for an archetypal physiological model, the Ultradian glucose-insulin model.
This observation suggests that DIU may profoundly affect clinical medical care, including glycemic management in the intensive care unit.
DIU may be relevant throughout biomedicine because delay is ubiquitous in physiological systems.
Developing DIU detection methods and assessing the impact of DIU on data assimilation techniques will be important future research directions.
Our work poses new mathematical questions at the interface of ergodic theory and infinite-dimensional delay dynamical systems.
}

\section{Introduction}
\label{s:intro}

Delay can significantly impact the dynamics of physiological systems at multiple scales.
At the level of genetic regulatory networks, distributed delay on the order of minutes results from the transcriptional, translational, and post-translational steps that lead to the production of functional regulator proteins.
Such delay can accelerate signaling within feedforward architectures~\cite{Josic2011stochastic}, alter the statistics of noise-induced switching phenomena~\cite{Gupta2013transcriptional, Kyrychko2018enhancing}, and produce oscillations in synthetic genetic circuits~\cite{Stricker2008fast}.
This paper is about a novel route through which delay can cause sustained temporal chaos within concrete dynamical systems of interest in physiology and biomedicine.
We call the resulting chaos \emph{delay-induced uncertainty} (DIU).
We believe that DIU has profound implications for oscillations that arise in systems physiology, especially the ultradian glucose-insulin oscillation observed within human endocrine physiology.

Clinical and laboratory practice throughout biochemistry, physiology, and medicine proceed from the assumption that the dynamics of measured quantities are not chaotic.
For instance, a clinician administers medication to a patient based on the belief that the medical intervention will not induce an unexpectedly erratic response.
While chaos has been observed in some physiological models~\cite{Abarbanel1993analysis, Li2004period, Mackey1977oscillation, Glass1988time, Glass1990chaos}, DIU in physiology is not yet well-understood by mathematicians, nor is its significance known to clinicians.
It is vital to uncover the mechanisms that produce DIU and develop precise mathematical characterizations of the resulting dynamics.
Moreover, it is crucial to assess the impact of DIU on data assimilation and clinical practice.
In this work, we uncover a route to DIU for oscillatory dynamical systems.

We formulate a general framework for the emergence of DIU in damped, driven oscillatory systems and then focus on the ultradian glucose-insulin oscillation.
The framework consists of three components.
First, delay renders the unforced dynamical system excitable.   
For the damped, driven oscillators we consider in this paper, delay produces a weakly stable limit cycle.
This can happen, for instance, via a delay-induced supercritical Hopf bifurcation~\cite{Karamched2015delayed, Song2005local, Wei2005hopf, Yan2006hopf, Bressloff2015frequency, Xu2017theory}.
Second, the unforced system possesses intrinsic shear.
For damped, driven oscillators, shear quantifies velocity gradients near the limit cycle.
Third, the forcing drive interacts with the shear to stretch and fold the phase space.
This interaction creates hyperbolicity in the dynamics and produces sustained temporal chaos.
The forcing drive does not overwhelm the intrinsic dynamics.
On the contrary, it interacts subtly with intrinsic shear to produce DIU.

We perform a number of experiments that show DIU emergence for an archetypal physiological model, the Ultradian glucose-insulin model~\cite{Sturis1991computer, Drozdov1995model, Keener1998mathematical, Albers2012population}.
This finite-dimensional model has been constructed to explain ultradian oscillations using a minimal number of components.
It includes compartments for insterstitial and plasma insulin, one glucose compartment, several feedback mechanisms that represent insulin-mediated glucose regulation by the pancreas, and hepatic responses.
Delayed regulatory feedback between insulin secretion and glucose released by the liver produces the oscillation.
The Ultradian model has been used to accurately predict glucose dynamics in humans~\cite{Albers2017personalized}.
It therefore provides an ideal setting for the investigation of DIU.

The presence of DIU in glucose-insulin dynamics may have profound implications for clinical care in the intensive care unit (ICU), where glucose and insulin treatments (external forcing drives) are central to glycemic management.
More generally, DIU is potentially relevant for any physiological system wherein delayed regulatory feedback controls try to maintain healthy homeostasis.
Examples include pulmonary and respiratory dynamics~\cite{Mackey1977oscillation, Sottile2018association}, cardiac dynamics~\cite{Christini2002introduction}, female endocrine dynamics~\cite{Graham2020reduced, Urteaga2019multi}, and neurological dynamics~\cite{Stroh2020estimating, Claassen2013nonconvulsive, Hodgkin1952quantitative}, to name but a few.
Indeed, the use of mathematical physiology within medicine has broad potential~\cite{Albers2018mechanistic, Zenker2007inverse}.

The Ultradian model is finite-dimensional because delayed regulatory feedback appears in the form of a three-stage filter.
Consequently, we use the theory of nonuniformly hyperbolic dynamical systems and specifically the theory of rank one maps~\cite{Wang2001strange, Wang2008toward, Wang2013dynamical} to precisely characterize DIU in the Ultradian system.
No theory of rank one maps for infinite-dimensional delay systems currently exists.
This suggests important open questions:
Do infinite-dimensional delay systems (delay differential equations) produce DIU?
How do we rigorously characterize DIU in this context?
We begin to answer the first of these questions here by showing numerically that DIU emerges in a delay variant of the linear shear flow model first studied by Zaslavsky~\cite{Zaslavsky1978simplest} and then by Lin and Young~\cite{Lin2008shear}.

We conclude the paper by discussing open mathematical questions inspired by DIU.
Further, we assess the potential impact of DIU on biomedicine and clinical practice.

\section{The Ultradian model}
\label{s:ultradian}

\begin{figure}[ht]
\centering
\includegraphics[width = 8cm]{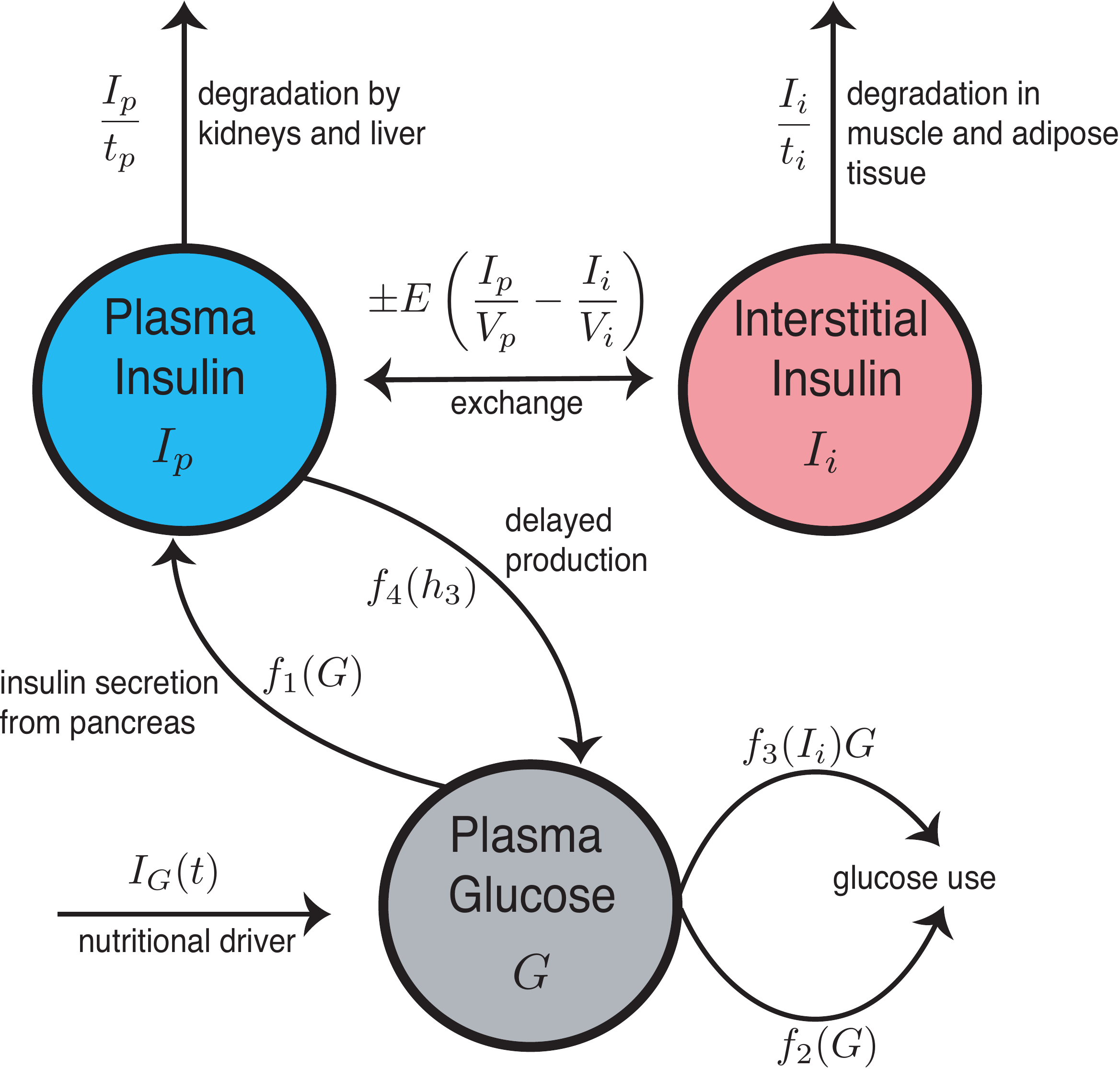}
\caption{Schematic of the Ultradian model of glucose-insulin dynamics.}
\label{high_level}
\end{figure}

We describe the Ultradian glucose-insulin model~\cite{Sturis1991computer, Drozdov1995model, Keener1998mathematical} as well as the external forcing drives that produce DIU. 

The Ultradian model is a compartment model with three state variables: plasma glucose ($G$), plasma insulin ($I_p$), and interstitial insulin ($I_i$).
See Fig.~\ref{high_level} for the model schematic. 
These three state variables are coupled to a three-stage linear delay filter, producing a $6$-dimensional phase space. 
The Ultradian model is particularly popular because it is the simplest physiological model that captures the main features of glucose-insulin oscillations~\cite{Sturis1991computer, Drozdov1995model} and provides a mechanistic description of the cause of the oscillations. 
The model includes two major negative feedback loops describing effects of insulin on glucose use and glucose production.
Both loops include glucose-based stimulation of insulin secretion. 

Oscillations in the Ultradian system depend on (i) a time delay of 30-45 minutes for the effect of insulin on glucose production and (ii) the slow effect of insulin on glucose use arising from insulin being in two distinct compartments. 
We focus on the former in this paper. 
Note that the Ultradian model includes physiological delay, but the system is \emph{finite-dimensional} because the delay assumes the form of a three-stage linear filter.

The full model is given by

\begin{subequations}
\label{eq:model1}
\begin{align}
\frac{\mathrm{d} I_p}{\mathrm{d} t} &= f_1(G) - E \left( \frac{I_{p}}{V_{p}}-\frac{I_i}{V_{i}} \right) - \frac{I_{p}}{t_{p}}
\\
\frac{\mathrm{d} I_i}{\mathrm{d} t} &= E \left( \frac{I_{p}}{V_{p}}-\frac{I_i}{V_{i}} \right) - \frac{I_{i}}{t_{i}}
\\
\frac{\mathrm{d} G}{\mathrm{d} t} &= f_4(h_3)+I_{G}(t)-f_2(G)-f_3(I_i)G
\\
\frac{\mathrm{d} h_1}{\mathrm{d} t} &= \frac{1}{t_d} ( I_p-h_1 )
\\
\frac{\mathrm{d} h_2}{\mathrm{d} t} &= \frac{1}{t_d} ( h_1-h_2 )
\\
\frac{\mathrm{d} h_3}{\mathrm{d} t} &= \frac{1}{t_d} ( h_2-h_3 ),
\end{align}
\end{subequations}
where $f_1(G)$ represents the rate of insulin production, $f_2(G)$ represents insulin-independent glucose use, $f_3(I_i)G$ represents insulin-dependent glucose use, and $f_4(h_3)$ represents delayed insulin-dependent glucose use.  The functional forms of these are

\begin{align*}
f_1(G) &= \frac{R_m}{1+ \exp \big( \frac{-G}{V_g c_1} + a_1 \big)} \\
f_2(G) &= U_b \Big( 1-\exp \Big( \frac{-G}{C_2 V_g} \Big) \Big) \\
f_3(I_i) &= \frac{1}{C_3 V_g} \Big( U_0 + \frac{U_m - U_0}{1 + (\kappa I_i)^{-\beta}} \Big) \\
f_4(h_3) &= \frac{R_g}{1 + \exp \big( \alpha \big( \frac{h_3}{C_5 V_p}  - 1 \big) \big) },
\end{align*}
with
\begin{align*}
\kappa &= \frac{1}{C_4} \left( \frac{1}{V_i} - \frac{1}{E t_i} \right).
\end{align*}
See Table~\ref{table:model_parameters} for parameter descriptions and nominal parameter values.

For our DIU experiments, we first consider an idealized nutritional driver $I_{G} (t)$ that includes a basal signal and pulsatile kicks.
The idealized nutritional driver is given by
\begin{equation}
\label{driver}
I_G(t) = I_0 + \sum_{n=0}^{\infty} A_n \delta(t - T_n),
\end{equation}
where $I_0$ is a basal nutritional input for the system, $T_n$ is the time of the $n$th feeding, and $A_n$ is the amount of carbohydrate in that meal. 
The signal $I_{G} (t)$ represents the external forcing in the Ultradian model.
The form of $I_{G} (t)$ in Eq.~\eqref{driver} produces the following dynamics:
Between two consecutive kicks ($T_{n-1} < t < T_{n}$), Ultradian dynamics evolve according to system~\eqref{eq:model1} with $I_{G} (t) = I_{0}$.
At the time $T_{n}$ of meal $n$, the glucose state variable, $G$, undergoes the instantaneous change $G \mapsto G + A_{n}$.
We demonstrate the emergence of DIU for both fixed and random kick amplitudes $(A_{n})_{n=0}^{\infty}$ and inter-kick times $(T_{n+1} - T_{n})_{n=0}^{\infty}$.

In reality, meals produce glucose kicks that are temporally localized but not instantaneous.
Our DIU results for the idealized nutritional driver strongly predict DIU emergence for the complex external forcing drives encountered in the intensive care unit.
To support this claim, we show that DIU remains present when we replace the $\delta $-kicks in Eq.~\eqref{driver} with square pulses of duration 30 minutes that arrive at 8 am, noon, and 6 pm.

\section{Delay-induced uncertainty for the Ultradian model}
\label{s:simulations}

We first formulate a general DIU framework for damped, driven oscillatory systems.
This framework explains the origins of sustained temporal chaos for the Ultradian model.
We then present a series of numerical experiments that demonstrate the presence and robustness of DIU in the Ultradian model.
We conclude Section~\ref{s:simulations} by providing a full dynamical profile of DIU in the Ultradian model and use ideas from smooth ergodic theory to support our numerical findings.

\subsection{DIU framework}

Our route to DIU for damped, driven oscillatory systems involves the following three components.

\begin{enumerate}[leftmargin=*, labelindent=\parindent, label=\textbf{(U\arabic*)}, ref=U\arabic*, topsep=1.5ex, itemsep=1ex]
\item
\label{li:excitability}
\textbf{Delay-induced excitability.}
Delay renders the unforced (intrinsic) dynamical system excitable by producing a weakly stable limit cycle.
As we will see, delay in the Ultradian model produces a limit cycle via a supercritical Hopf bifurcation.

\item
\label{li:intrinsic-shear}
\textbf{Intrinsic shear.}
Shear refers to significant velocity gradients in a tubular neighborhood of the limit cycle.
Atmospheric wind shear provides a good mental picture of the phenomenon.

\item
\label{li:external-forcing}
\textbf{External forcing allows shear to act.}
External forcing allows the shear to stretch and fold the phase space, thereby creating hyperbolicity in the dynamics.

\end{enumerate}

\begin{figure}[ht]
\centering
\includegraphics[width = 15cm]{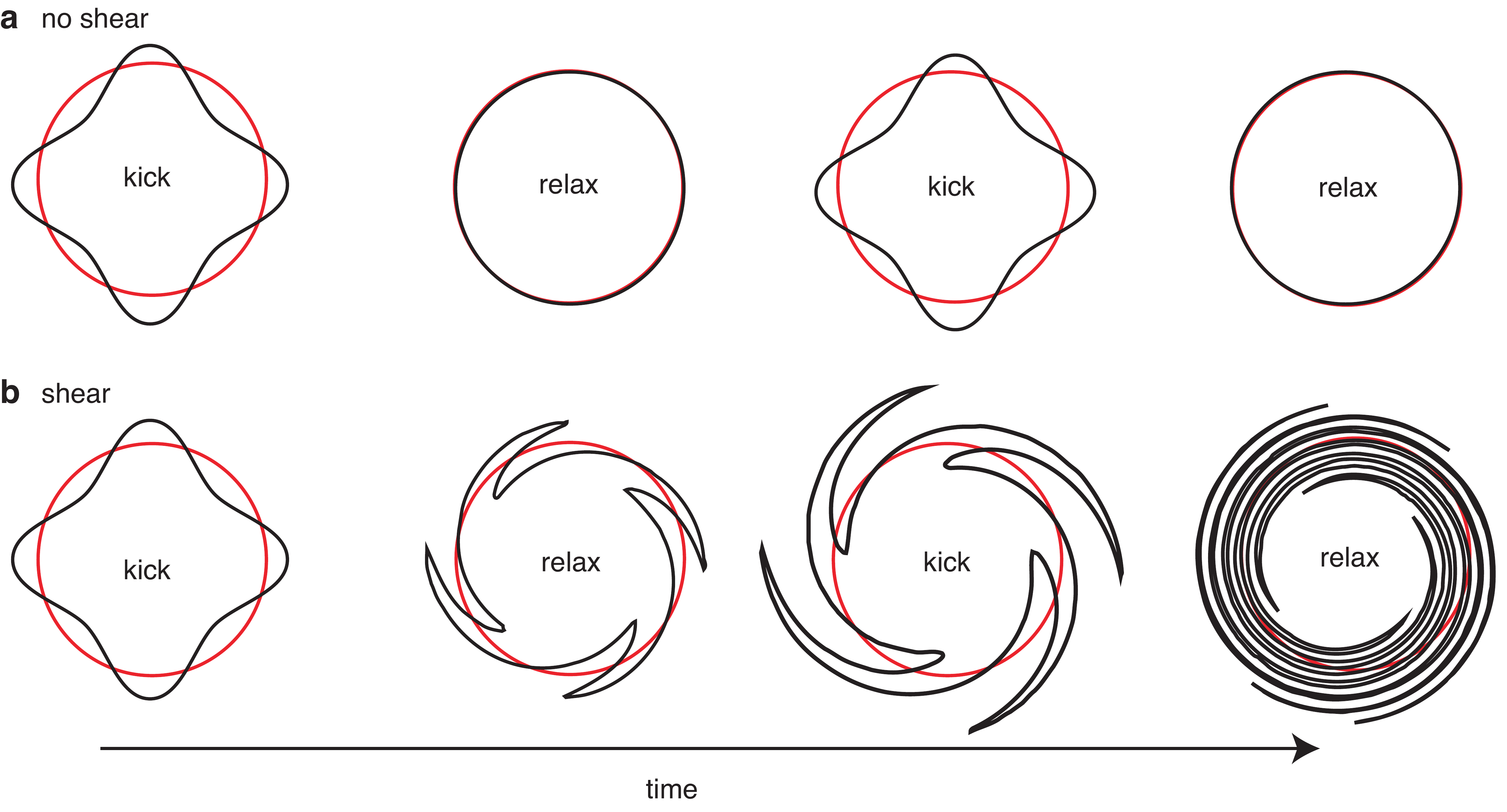}
\caption{\textbf{Kick-relaxation dynamics of the Ultradian model with pulsatile forcing.}
When the system is kicked, the limit cycle (red circle) deforms (black curve).
Before the next kick, the system relaxes toward the limit cycle.
\textbf{(a)}
In the absence of shear, the kicked limit cycle quickly relaxes.
The phase space does not stretch and fold.
\textbf{(b)}
When shear is strong, the kicked limit cycle stretches and folds during the relaxation phase (assuming the time between kicks is long enough to allow stretch and fold geometry to manifest).
The kick-relaxation cycle produces sustained temporal chaos.
}
\label{fig:sawtooth}
\end{figure}

Figure~\ref{fig:sawtooth} illustrates the geometric mechanism behind the emergence of sustained temporal chaos.
Since our glucose forcing signals~\eqref{driver} are pulsatile, the evolution of the Ultradian system decomposes into windows of relaxation punctuated by kicks.
The amount of shear near the limit cycle determines how the kick-relaxation cycle acts on phase space.
In the absence of shear (Figure~\ref{fig:sawtooth}a), the kicked limit cycle would calmly relax before the next kick.
The phase space would not stretch and fold in this case.
When shear is strong, the phase space would stretch and fold between kicks (Figure~\ref{fig:sawtooth}b), producing sustained temporal chaos.

Importantly, DIU is not a phenomenon wherein the external forcing simply overwhelms the intrinsic dynamics.
On the contrary, forcing amplitudes can be quite small.
Forcing acts as an amplifier in our DIU framework, amplifying the impact of intrinsic shear to produce rich, complex dynamics.

\subsection{Simulation results for idealized pulsatile forcing}

We deploy our DIU framework to establish DIU emergence in the Ultradian model with glucose input signal~\eqref{driver}.
We show that DIU emerges for both constant and random kick amplitudes and inter-kick times.
The results of this section are numerical in nature for the following reason.
Since the Ultradian model is finite-dimensional, the analysis of kicked limit cycles for nonlinear systems in~\cite{Ott2010from, Wang2003strange} provides rigorous mathematical support for our simulation results.
This analysis leverages the theory of rank one maps~\cite{Wang2001strange, Wang2008toward, Wang2013dynamical} to provide a dynamical profile of DIU.
This profile includes the existence of a strange attractor.
However, the theory of rank one maps has not yet been extended to treat random inter-kick times or random kick amplitudes.
Even when both kick amplitude and inter-kick time are held constant ($A_{n} = A$ and $T_{n} = nT$ for all $n$), the current theory of rank one maps cannot tell us if a strange attractor exists for \emph{specific} values of $A$ and $T$~\cite{Ott2010from, Wang2003strange}.
Rather, the rigorous applications of the theory developed thus far prove the existence of strange attractors for parameter sets of positive Lebesgue measure.
Moreover, the theory of~\cite{Ott2010from, Wang2003strange} is an asymptotic theory whereas our numerical experiments take place in a practical parameter regime.

We therefore analyze the Ultradian model numerically and use the maximal Lyapunov exponent, $\Lambda_{\max}$, as a DIU diagnostic:
$\Lambda_{\max} > 0$ indicates DIU, while $\Lambda_{\max} < 0$ indicates that DIU is absent.

\vspace{1.5ex}
\noindent
\textbf{Parameter selection.}
Excluding the nutritional driver, we set all Ultradian model parameters to the values in Table~\ref{table:model_parameters} for our simulations.
For the nutritional driver, we set the basal rate $I_{0}$ to zero.
We are therefore free only to tune the delay $t_{d}$ and choose models for the kick amplitudes $A_{n}$ and the inter-kick times $T_{n+1} - T_{n}$.

\vspace{1.5ex}
\noindent
\textbf{DIU framework for the Ultradian model.}
The emergence of a limit cycle in system~\eqref{eq:model1} as the delay $t_{d}$ increases invokes framework~(\ref{li:excitability})--(\ref{li:external-forcing}) for the presence of DIU in the Ultradian system.

\begin{figure}[ht]
\centering
\includegraphics[width = 15cm]{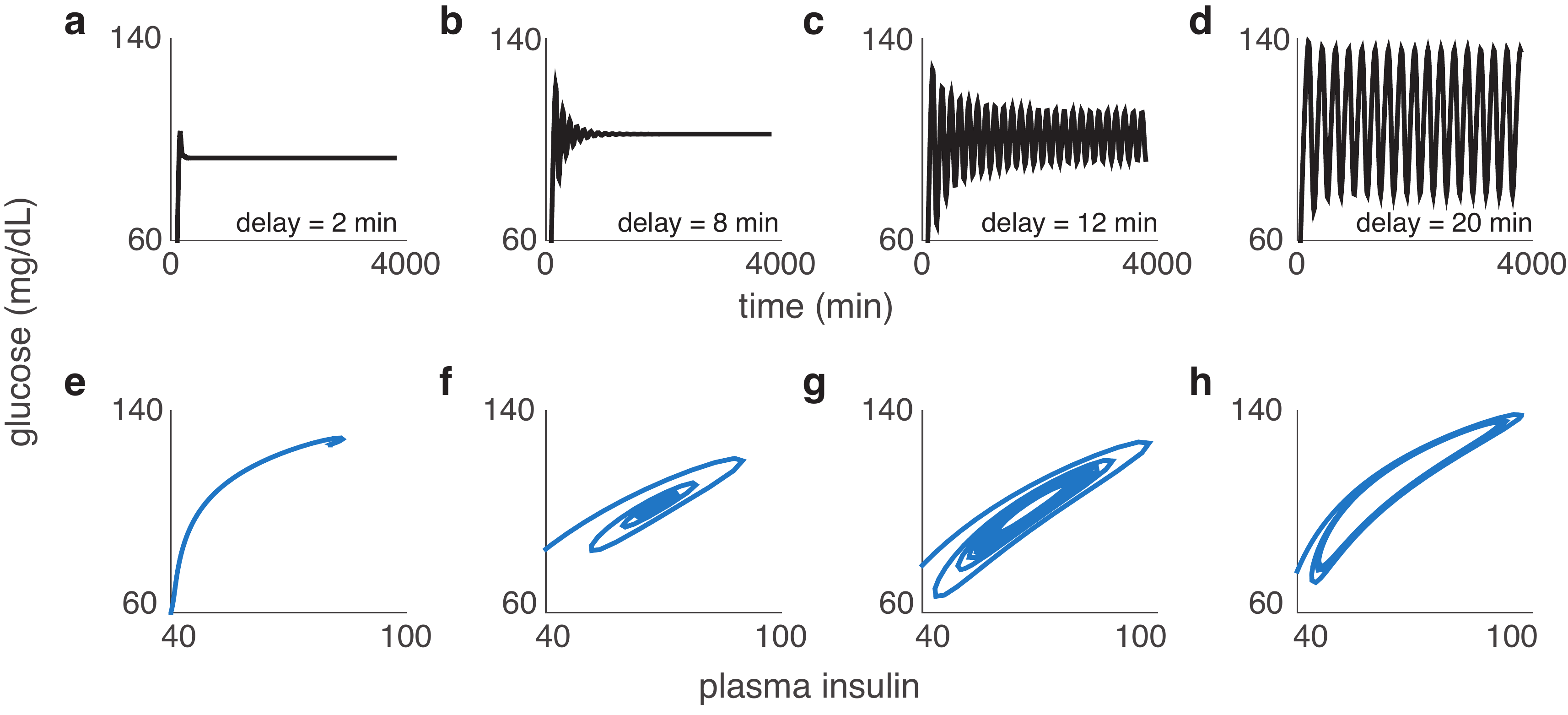}
\caption{
For the unforced Ultradian model, a supercritical Hopf bifurcation produces a limit cycle as delay $t_{d}$ increases. 
Top row: Glucose timeseries for \textbf{(a)} $t_d = 2$ min, \textbf{(b)} $t_d = 8$ min, \textbf{(c)} $t_d = 12$ min, and \textbf{(d)} $t_d = 20$ min. 
Bottom row: Projection of a phase space trajectory onto glucose-plasma insulin space  for \textbf{(e)} $t_d = 2$ min, \textbf{(f)} $t_d = 8$ min, \textbf{(g)} $t_d = 12$ min, and \textbf{(h)} $t_d = 20$ min. 
Other parameter values are given in Table~\ref{table:model_parameters}.
}
\label{fig:delay_hopf}
\end{figure}

(\ref{li:excitability}).
Consider the unforced version of system~\eqref{eq:model1}, obtained by removing $I_{G} (t)$.
For a variety of time delays, Figure~\ref{fig:delay_hopf} shows glucose timeseries (top row) and two-dimensional projections of phase space trajectories (bottom row) generated by the unforced Ultradian system.
We see that a stable equilibrium bifurcates into a limit cycle as the system undergoes a supercritical Hopf bifurcation at a delay value $t_{d}^{*}$ satisfying $8 < t_{d}^{*} < 12$.
The presence of the limit cycle implies excitability.

(\ref{li:intrinsic-shear}) and (\ref{li:external-forcing}).
We claim that limit cycles subjected to pulsatile forcing drives generically satisfy these two framework components.
The geometric ideas of Wang and Young~\cite{Wang2003strange} and the quantitative analysis of Ott and Stenlund~\cite{Ott2010from} support this claim for finite-dimensional nonlinear systems.
For~(\ref{li:intrinsic-shear}), shear can be understood geometrically by examining the shape of the strong stable foliation in a tubular neighborhood of the limit cycle~\cite{Wang2003strange}.
Ott and Stenlund~\cite{Ott2010from} quantify shear by defining a \emph{shear integral} that represents the accumulation of shear as one traverses the limit cycle.

For~(\ref{li:external-forcing}), the external forcing must interact with the shear in order to produce sustained temporal chaos.
We claim that this happens for generic pulsatile forcing drives.
To support this claim, Wang and Young prove that given a $C^{4}$ flow on a Riemannian manifold that admits a hyperbolic limit cycle, periodic kicks will produce strange attractors for an open set of $C^{3}$ kick functions (Theorem~1 of~\cite{Wang2003strange}).
Ott and Stenlund~\cite{Ott2010from} define a function that quantifies the forcing-shear interaction and assume that this function is Morse in their main theorem on the existence of strange attractors.
They conjecture that this assumption will hold for a generic kick-generating vector field, both in terms of topological genericity and prevalence.
See Remark~2.1 of~\cite{Ott2010from} for a discussion of the conjecture and~\cite{Ott2005prevalence} for information about prevalence.
Note that for the Ultradian model, the kicks provided by the nutritional driver have no spatial variation \emph{with respect to the original state variables}.
Such variation should be present, however, \emph{for the coordinate system near the limit cycle} developed in~\cite{Ott2010from}.

\vspace{1.5ex}
\noindent
\textbf{The maximal Lyapunov exponent as a diagnostic tool.}
We now present the numerical experiments that establish the presence of DIU for the Ultradian model.
We use the maximal Lyapunov exponent as a DIU diagnostic: $\Lambda_{\max} > 0$ indicates DIU, while $\Lambda_{\max} < 0$ indicates its absence.
We compute the maximal Lyapunov exponent by solving system~\eqref{eq:model1} as follows.
During the relaxation intervals $(T_{n-1}, T_{n})$ between kicks, we integrate the differential equations using MATLAB's ode23s stiff solver.
At kick times $T_{n}$, we pause the differential equation solver and apply the diffeomorphism of phase space induced by the kick $G \mapsto G + A_{n}$.
We compute $\Lambda_{\max}$ by completing $10^{5}$ kick-relaxation cycles.
Our maximal Lyapunov exponent therefore quantifies the amount of expansion per kick-relaxation cycle.

\vspace{1.5ex}
\noindent
\textbf{Constant kick amplitude, periodic or Poissonian kicks.}
For our first set of experiments, we choose a value of the delay $t_{d}$ such that the limit cycle is present in the unforced system, and then hold $t_{d}$ fixed.
We consider kicks of constant intensity, $A_{n} = A$ for all $n$.
Kick times are either periodic, $T_{n} = nT$ for all $n$, or Poissonian.
In the Poissonian case, the inter-kick times $T_{n+1} - T_{n}$ are independent and exponentially distributed with mean $T$.
We show that DIU emerges even for these relatively simple forms of the nutritional driver~\eqref{driver} by examining how the maximal Lyapunov exponent depends on $A$ and $T$.

\begin{figure}[ht]
\centering
\includegraphics[width = 15cm]{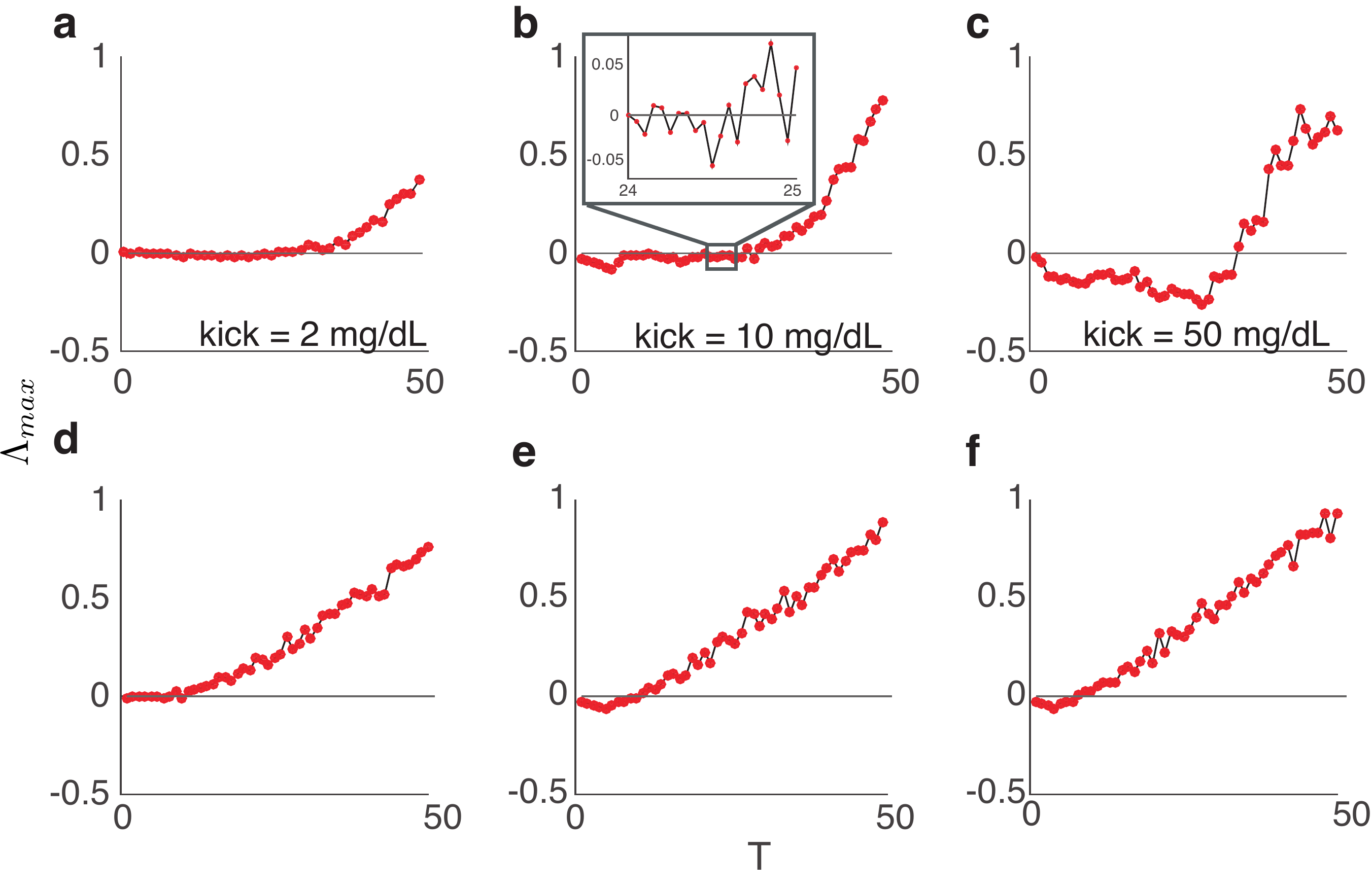}
\caption{
\textbf{Maximal Lyapunov exponent in the Ultradian model as a function of kick timing.}
Positive $\Lambda_{\max}$ indicates the presence of DIU.
\textbf{(a)-(c)}
For periodic kicks, plots of $\Lambda_{\max}$ versus the inter-kick time $T$ for several fixed values of kick amplitude $A$.
\textbf{(a)}
$A = 2$ mg/dL.
\textbf{(b)}
$A =10$ mg/dL.
\textbf{(c)}
$A = 50$ mg/dL.
\textbf{(d)-(f)}
For Poissonian kicks, plots of $\Lambda_{\max}$ versus mean inter-kick time $T$ for several fixed values of kick amplitude $A$.
\textbf{(d)}
$A = 2$ mg/dL.
\textbf{(e)}
$A =10$ mg/dL.
\textbf{(f)}
$A = 50$ mg/dL.
Other parameter values are as in Table~\ref{table:model_parameters}.
}
\label{fig:lyapounoff1}
\end{figure}

\begin{figure}[ht]
\centering
\includegraphics[width = 15cm]{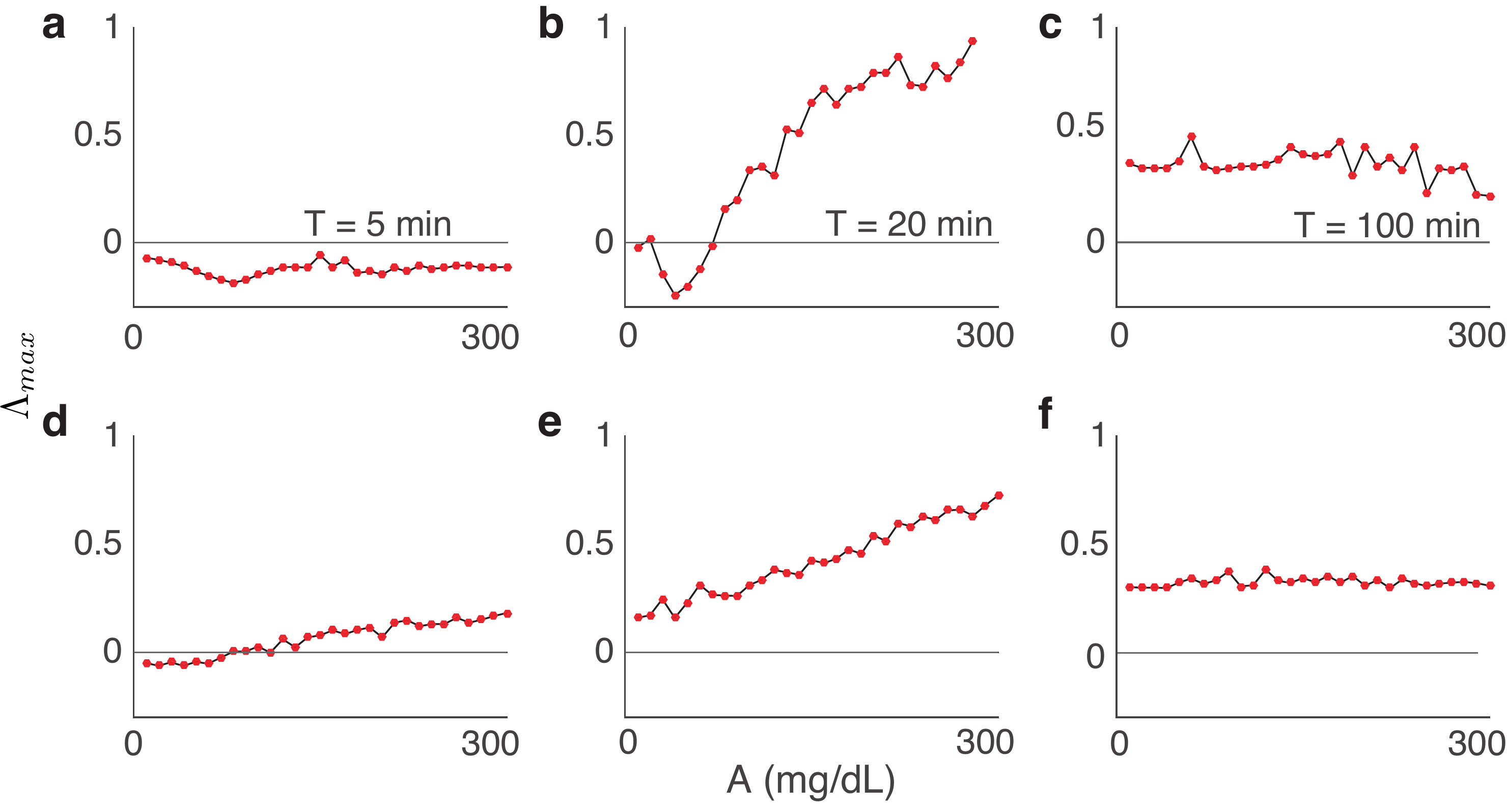}
\caption{
\textbf{Maximal Lyapunov exponent in the Ultradian model as a function of kick amplitude.}
Positive $\Lambda_{\max}$ indicates the presence of DIU.
\textbf{(a)-(c)}
For periodic kicks, plots of $\Lambda_{\max}$ as a function of kick amplitude for several different values of inter-kick time $T$.
\textbf{(a)}
$T = 5$ min.
\textbf{(b)}
$T = 20$ min.
\textbf{(c)}
$T = 100$ min.
\textbf{(d)-(f)}
For Poissonian kicks, plots of $\Lambda_{\max}$ as a function of kick amplitude for several different values of mean inter-kick time $T$.
\textbf{(d)}
$T = 5$ min.
\textbf{(e)}
$T = 20$ min.
\textbf{(f)}
$T = 100$ min.
Other parameter values are as in Table~\ref{table:model_parameters}.}
\label{fig:lyapounoff2}
\end{figure}

We compute the maximal Lyapunov exponent as follows.
For simulations involving periodic kicks (see Figs.~\ref{fig:lyapounoff1}a-c, \ref{fig:lyapounoff2}a-c), we track two solutions to system~\eqref{eq:model1}, initially separated by $d_0 = 10^{-8}$.
Think of one of these solutions as the base solution and the other as a secondary, perturbed solution.
After the first kick-relaxation cycle, we compute the separation $d_1$ between the solutions at time $T$ and store the quantity $\log(d_1/d_0)$ in a vector. 
We then renormalize by moving the secondary orbit toward the base orbit so that the distance between the two resets to $d_{0}$.
We proceed in this manner for $10^{5}$ kick-relaxation cycles. 
This produces a vector containing $10^{5}$ values of $\log (d_1/d_0)$.
Averaging over this vector produces $\Lambda_{\max}$, the maximal Lyapunov exponent.

For simulations involving Poissonian kicks (see Figs.~\ref{fig:lyapounoff1}d-f, \ref{fig:lyapounoff2}d-f), the maximal Lyapunov exponent is a random variable, as it depends \emph{a priori} on the random inter-kick times.
To compute it, we first sample $10^5$ inter-kick times from the exponential distribution with mean $T$.
These samples produce a single realization of the stochastic process.
We compute the maximal Lyapunov exponent for this realization by proceeding as we did in the case of periodic kicks.
That is, we compute $\log(d_1/d_0)$ following each kick-relaxation cycle and then average.
Finally, we average the realization-dependent maximal Lyapunov exponent over $10^{5}$ realizations of the Poisson process.
Abusing notation slightly, we call this average $\Lambda_{\max}$.

Figures~\ref{fig:lyapounoff1}a-c and \ref{fig:lyapounoff2}a-c display maximal Lyapunov exponent results for the case of constant kick amplitude and periodic kicks.
Here, $\Lambda_{\max}$ is a function of the kick amplitude $A$ and the inter-kick time $T$.

For three different fixed values of $A$, $\Lambda_{\max}$ becomes positive as $T$ increases, indicating the onset of DIU (Figure~\ref{fig:lyapounoff1}a-c).
This is consistent with the intuition from Figure~\ref{fig:sawtooth}: Larger values of $T$ allow more time for the phase space to stretch and fold between kicks.
The maximal Lyapunov exponent depends on $T$ in a particularly interesting way when $A = 50$ mg/dL (Figure~\ref{fig:lyapounoff1}c).
Here, $\Lambda_{\max} < 0$ for small values of $T$, indicating that DIU is absent and suggesting that the time-$T$ map of the system possesses an attractor that is diffeomorphic to the limit cycle present in the unforced Ultradian system.
By contrast, $\Lambda_{\max} > 0$ for large values of $T$, indicating the presence of DIU and suggesting that the time-$T$ map of the system possesses a strange attractor.
The inset in Figure~\ref{fig:lyapounoff1}b shows that $\Lambda_{\max}$ fluctuates around zero for moderately large values of $T$.
This suggests that the time-$T$ map of the system possesses horseshoes (transient chaos).

In Figure~\ref{fig:lyapounoff2}a-c, we compute $\Lambda_{\max}$ as a function of $A$ for different fixed values of $T$.
When $T$ is small ($T = 5$ min), the maximal Lyapunov exponent is negative for all of the values of $A$ we have simulated, indicating that DIU is absent, robustly so with respect to $A$, when $T$ is small (Figure~\ref{fig:lyapounoff2}a).
By contrast, when $T$ is large ($T = 100$ min), $\Lambda_{\max}$ is positive for all of the values of $A$ we have tested, indicating that DIU is present even when $A$ is small (Figure~\ref{fig:lyapounoff2}c).
We observe a transition from quiescence to DIU as $A$ increases when $T$ is moderately large ($T = 20$ min, Figure~\ref{fig:lyapounoff2}b).

Figures~\ref{fig:lyapounoff1}d-f and \ref{fig:lyapounoff2}d-f display maximal Lyapunov exponent results for the case of constant kick amplitude and Poissonian kicks.
Here, $\Lambda_{\max} > 0$ for most of the kick amplitudes $A$ and mean inter-kick times $T$ we tested, indicating robust presence of DIU.
Note that when $A$ is fixed at $A = 50$ mg/dL, DIU onset occurs significantly earlier in the Poissonian case than in the periodic case as $T$ increases (Figure~\ref{fig:lyapounoff1}c,f).

\begin{figure}[ht]
\centering
\includegraphics[width = 0.95\textwidth]{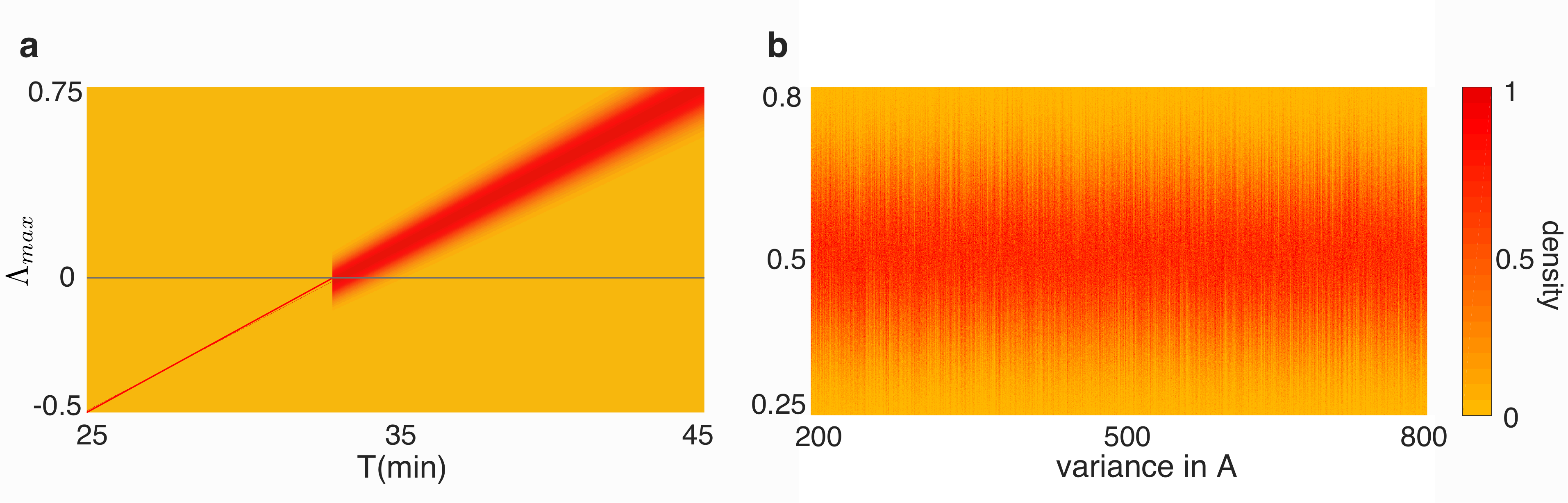}
\caption{
\textbf{The distribution of} $\bm{\Lambda_{\max}}$ \textbf{for the Ultradian model when the kicks are independent and uniformly distributed.}
Once again we see evidence of DIU.
\textbf{(a)}
We fix the kick amplitude distribution $A \sim U[45, 55]$, so that $\langle A \rangle = 50$ mg/dL and $\Var [A] = 8.\bar{3}$.
For $T$ values that produce $\Lambda_{\max} < 0$ in Figure~\ref{fig:lyapounoff1}c, the distribution of $\Lambda_{\max}$ resembles a Dirac delta.
However, for $T$ values that produce $\Lambda_{\max} > 0$ in Figure~\ref{fig:lyapounoff1}c, a distribution of values emerges for $\Lambda_{\max}$ that broadens as $T$ grows.  
\textbf{(b)}
We fix $T = 40$ min, fix $\langle A\rangle = 50$ mg/dL, and plot the distribution of $\Lambda_{\max}$ as a function of $\Var [A]$.
}
\label{fig:random_a}
\end{figure}

\vspace{1.5ex}
\noindent
\textbf{Uniformly distributed kick amplitudes, periodic kicks.}
We claim that DIU emergence is robust - DIU will emerge regardless of the particular shape of the pulsatile forcing.
Our remaining experiments with the Ultradian model support this claim.
For the next set of experiments, we make the (more realistic) assumption that kick amplitudes are random, rather than constant.
In particular, we assume the kick amplitudes $A_{n}$ are independent and uniformly distributed, while the kicks are periodic in time with inter-kick time $T$.
Figure~\ref{fig:random_a}a shows the distribution of $\Lambda_{\max}$ as a function of $T$ when the kick amplitudes are drawn from the uniform distribution on $[45,55]$.
When $T$ is small, the distribution of $\Lambda_{\max}$ is essentially a Dirac delta at a negative value.
Interestingly, at the moment $E[\Lambda_{\max}]$ crosses zero, the variance of $\Lambda_{\max}$ immediately becomes positive, and continues to grow as $T$ increases.

For our second experiment with uniformly distributed kick amplitudes, we fix $T$ at $T = 40$ min, fix the mean of the kick amplitude distribution at $50$ mg/dL, and examine how the distribution of $\Lambda_{\max}$ varies with the variance of the kick amplitude distribution.
Notice that $T = 40$ min is beyond the critical value at which we see an abrupt behavioral change in Figure~\ref{fig:random_a}a.
Interestingly, the overall width of the $\Lambda_{\max}$ distribution seems to be insensitive to kick amplitude variance, yet we see subtle variation at fine scales (Figure~\ref{fig:random_a}b).

\begin{figure}[ht]
\centering
\includegraphics[width=0.95\textwidth]{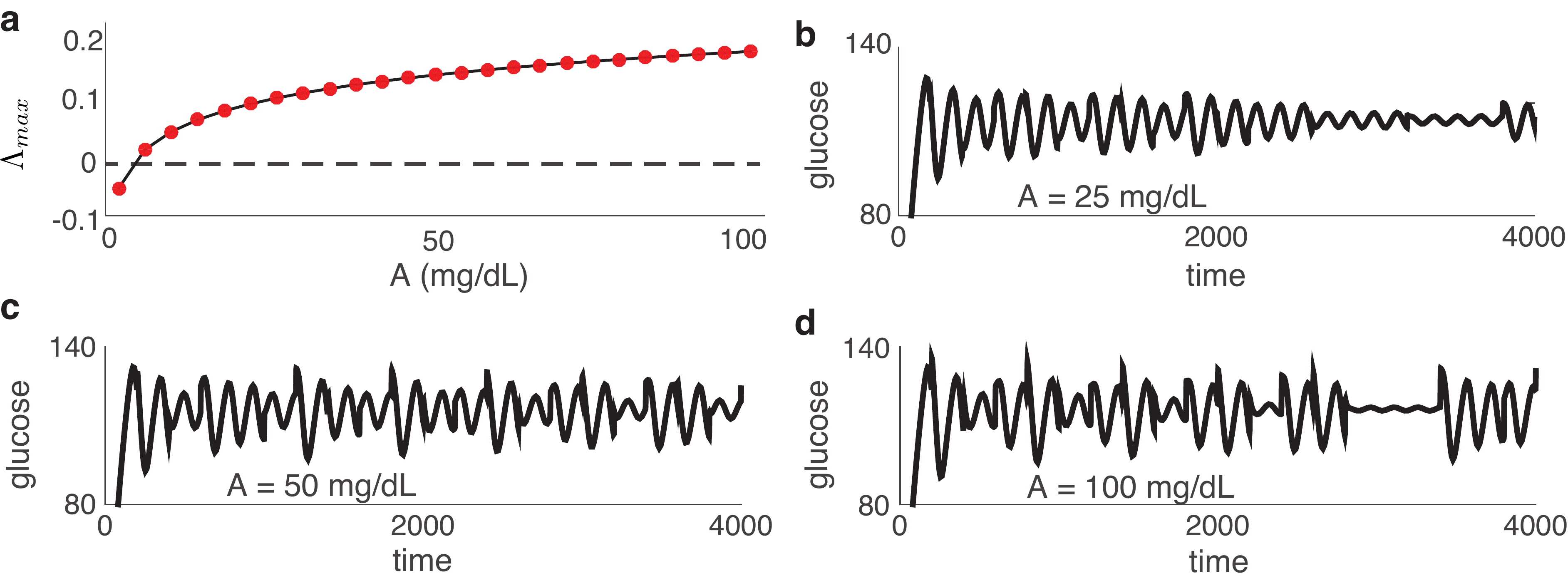}
\caption{
\textbf{Robust DIU emergence for a meal-like carbohydrate signal.}
We replace the sum of $\delta$-pulses in the nutritional driver~\eqref{driver} with square pulses of height $A$ that have a duration of 30 minutes and begin daily at 8 am, noon, and 6 pm.
\textbf{(a)}
The top Lyapunov exponent is positive even for small values of $A$, indicating robust emergence of DIU for this forcing signal.
\textbf{(b)-(d)}
Sample glucose timeseries for three values of $A$ in the DIU regime.
}
\label{fig:square-pulses}
\end{figure}

\vspace{1.5ex}
\noindent
\textbf{A meal-like carbohydrate input signal.}
For our final experiment with the Ultradian model, we replace the sum of $\delta$-pulses in the nutritional driver~\eqref{driver} with square pulses of height $A$ that have a duration of $30$ minutes and begin daily at $8$ am, noon, and $6$ pm.
Figure~\ref{fig:square-pulses}a shows that DIU quickly emerges as $A$ increases.
For values of $A$ in the DIU regime, the corresponding glucose timeseries behave in an interesting way (Figure~\ref{fig:square-pulses}b-d).
When $A = 100$, for instance, the glucose timeseries contains windows with erratic behavior and windows wherein the glucose signal is nearly constant.
Such behavior would potentially confuse clinicians and can complicate data-based detection and interpretation of DIU.

\subsection{Dynamical profile of DIU for the Ultradian model}

Since the Ultradian model is finite-dimensional, the analysis of kicked limit cycles for nonlinear systems in~\cite{Ott2010from, Wang2003strange} supports our numerical findings and leverages the theory of rank one maps~\cite{Wang2001strange, Wang2008toward, Wang2013dynamical} to provide a dynamical profile of DIU.
We associate DIU with the existence of a strange attractor that supports a unique ergodic Sinai-Ruelle-Bowen measure.
The system has a positive Lyapunov exponent (sustained temporal chaos) and possesses rich statistical properties.
These include a dynamical version of the central limit theorem, exponential decay of correlations, and a large deviation principle.

\begin{figure}[ht]
\centering
\includegraphics[width = 15cm]{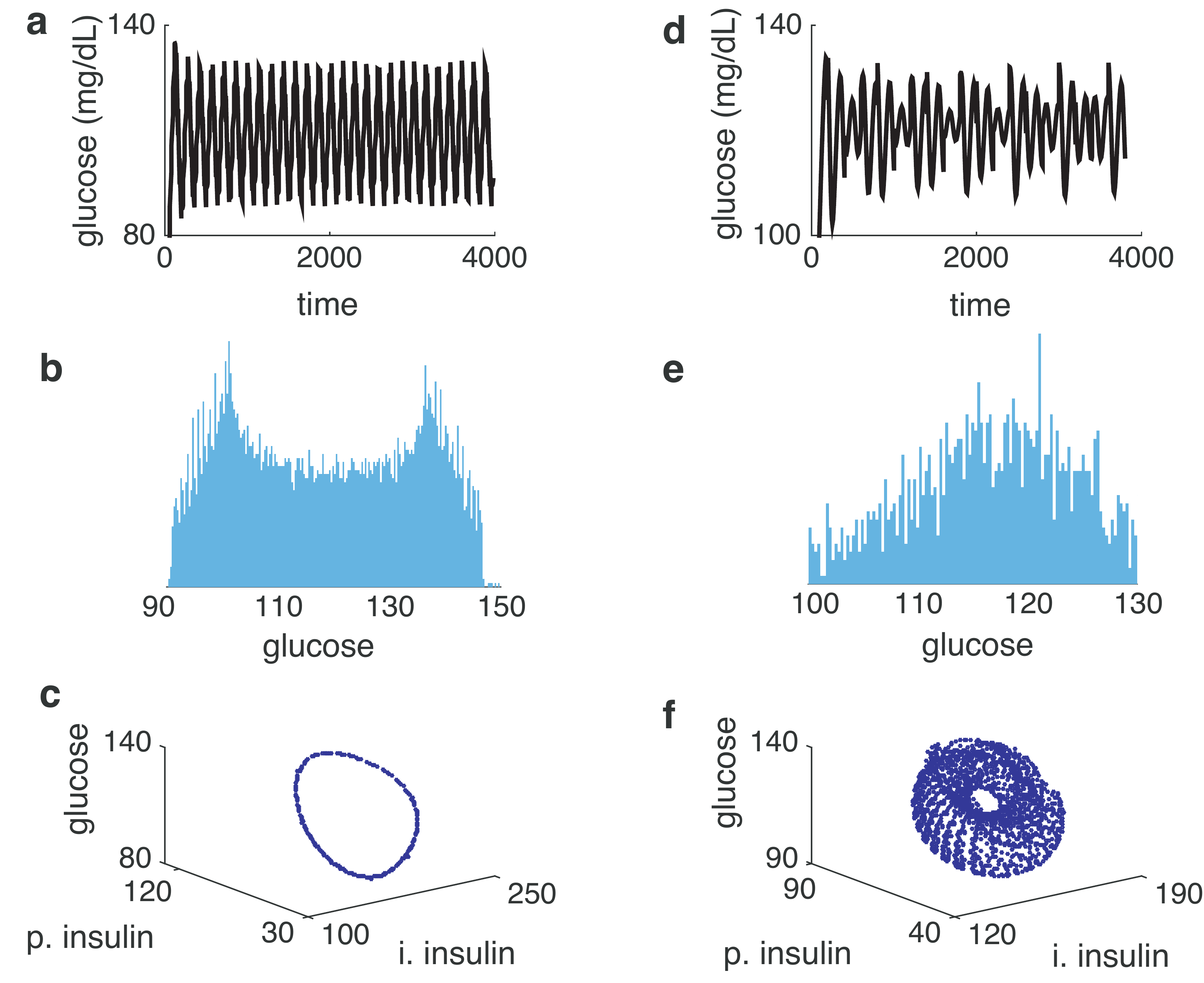}
\caption{
\textbf{Dynamical profile of Ultradian dynamics in the absence (left column) and presence (right column) of DIU.}
\textbf{(a)}
Glucose levels oscillate in a regular manner.
\textbf{(b)}
The empirical glucose distribution associated with the timeseries from~(a) is bimodal.
\textbf{(c)}
The time-$T$ map generated by the Ultradian system possesses an attractor that is diffeomorphic to the limit cycle of the unforced system.
\textbf{(d)}
Glucose levels evolve erratically.
\textbf{(e)}
The empirical glucose distribution associated with the timeseries from~(d) is unimodal and appears to be approximately Gaussian but with finite support.
\textbf{(f)}
The time-$T$ map generated by the Ultradian system possesses a strange attractor.
Left column: $T = 20$ min.
Right column: $T = 200$ min.
All panels: $t_d = 12$ min, $A = 10$ mg/dL, and other parameter values are as in Table~\ref{table:model_parameters}.
}
\label{fig:invariant_density}
\end{figure}

Figure~\ref{fig:invariant_density} illustrates how DIU impacts Ultradian dynamics.
Here, we simulate the Ultradian model with delay $t_{d} = 12$ min, a value for which the unforced system possesses a limit cycle.
We drive the system using nutritional driver~\eqref{driver} with periodic glucose kicks of constant amplitude $A = 10$ mg/dL.
We select a small value of the inter-kick time for which DIU is absent ($T = 20$ min, left column) and a larger value for which DIU is present ($T = 200$ min, right column).

The first two rows of Figure~\ref{fig:invariant_density} show representative glucose timeseries and corresponding empirical glucose distributions.
When DIU is absent, glucose levels oscillate regularly as expected, but interestingly the empirical glucose distribution is bimodal.
By contrast, the erratic behavior of the glucose timeseries in the presence of DIU reflects the chaos in the system.
Notice that the empirical distribution in the DIU case appears to be approximately Gaussian.
This observation is consistent with the results on Sinai-Ruelle-Bowen measures for general kicked limit cycles in~\cite{Ott2010from, Wang2003strange}, where scalar observables are shown to satisfy a dynamical version of the central limit theorem.

Figure~\ref{fig:invariant_density}d,e suggests that the impact of DIU on clinical practice will be subtle and complex.
When employing a \emph{single-orbit perspective} on dynamics, sustained temporal chaos renders rational medical intervention extremely difficult.
If only the \emph{statistical} behavior of observables of the dynamics (such as glucose level) matters in a particular setting, then DIU may be \emph{beneficial}, since the results of~\cite{Ott2010from, Wang2003strange} suggest that observables of Ultradian dynamics behave with a high level of statistical regularity.

We plot the attractors of the time-$T$ map generated by the Ultradian system in Figure~\ref{fig:invariant_density}c,f.
When DIU is absent, the attractor is diffeomorphic to the limit cycle of the unforced system.
In the presence of DIU, we observe a strange attractor with intricate geometry.
These results are consistent with the rigorous theory of~\cite{Ott2010from, Wang2003strange}.

\section{Delay-induced uncertainty for delay linear shear flow}
\label{s:linear-model}

There exists no theory of rank one maps for infinite-dimensional delay systems at this time.
Nevertheless, the geometric principles behind our DIU framework remain valid for infinite-dimensional dynamics.
We believe that it will be possible to develop a rigorous DIU theory for delay differential equations (DDEs).
Such a theory would have considerable value given the ubiquity of DDE modeling in mathematical physiology.

We show that it is possible for delay differential equations to produce DIU by demonstrating that even a simple DDE does so.

\subsection{Delay linear shear flow}
The dynamics take place on the cylinder $\mathbb{S}^{1} \times \mathbb{R}$.
Writing $\theta$ for the $\mathbb{S}^{1}$-coordinate and $z$ for the $\mathbb{R}$-coordinate, delay linear shear flow is generated by the delay differential equations
\begin{subequations}
\label{e:dlsf-kicked}
\begin{align}
\dot{z} (t) &= - \lambda z (t - \tau) + A \Phi (\theta) \sum_{n=0}^{\infty} \delta (t - nT),
\label{e:dlsf-kicked-z}
\\
\dot{\theta} &= 1 + \sigma z.
\label{e:dlsf-kicked-theta}
\end{align}
\end{subequations}
System~\eqref{e:dlsf-kicked} is infinite-dimensional due to the delay $\tau$: One must specify an initial history $h : [-\tau, 0] \to \mathbb{S}^{1} \times \mathbb{R}$ in order to propagate the dynamics forward.
The final term on the right side of Eq.~\eqref{e:dlsf-kicked-z} represents periodic pulsatile kicks.
Here $A \geqs 0$ is the kick amplitude, $T > 0$ is the inter-kick time, $\Phi : \mathbb{S}^{1} \to \mathbb{R}$ describes the kick profile, and $\delta$ is the Dirac delta.
The Dirac delta has the following interpretation:
At each nonnegative integer multiple of $T$, each point in $\mathbb{S}^{1} \times \mathbb{R}$ instantaneously moves from $(\theta, z)$ to $(\theta, z + A \Phi (\theta))$.
Between kicks, the dynamics are governed by the unforced delay equations
\begin{subequations}
\label{e:dlsf-intrinsic}
\begin{align}
\dot{z} (t) &= - \lambda z (t - \tau),
\\
\dot{\theta} &= 1 + \sigma z.
\label{e:dlsf-intrinsic-thetaDot}
\end{align}
\end{subequations}
Together with delay $\tau$, the parameters $\lambda > 0$ and $\sigma \in \mathbb{R}$ shape the dynamics of the unforced dynamical system.

\begin{figure}[ht]
\centering
\includegraphics[width = 15cm]{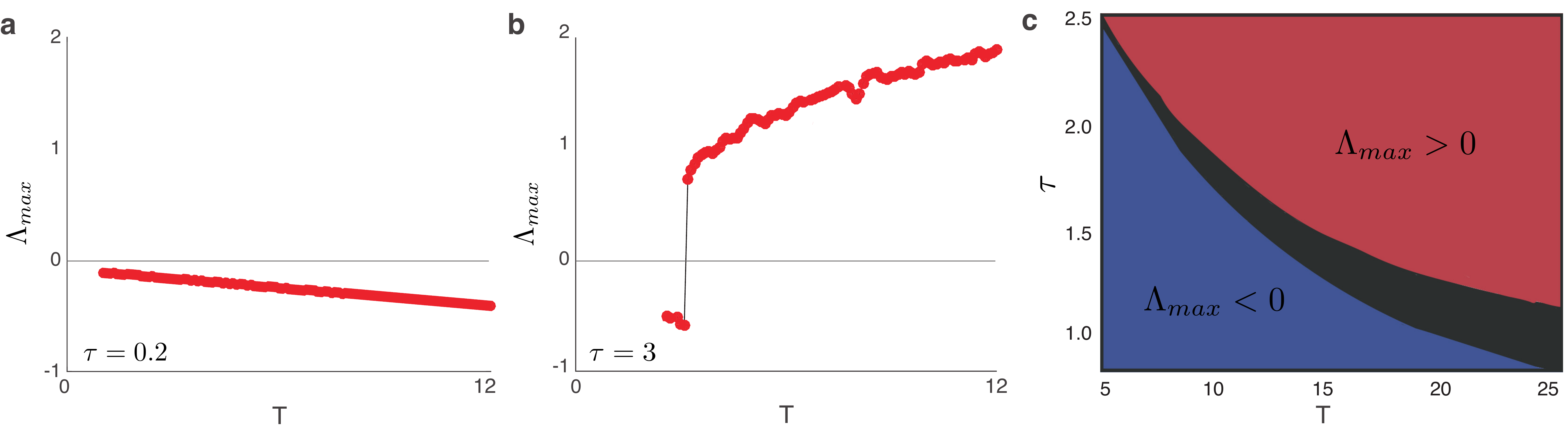}
\caption{
{\bfseries DIU for delay linear shear flow~\eqref{e:dlsf-kicked}.}
\textbf{(a)}
When delay $\tau$ is small, the top Lyapunov exponent $\Lambda_{\mathrm{max}}$ is negative over the tested range of inter-kick times $T$.
DIU has not emerged when delay is small.
\textbf{(b)}
When delay is substantial, DIU emerges when $\Lambda_{\mathrm{max}}$ transitions from negative to positive.
\textbf{(c)}
Heatmap illustrating the sign of $\Lambda_{\mathrm{max}}$ as a function of the delay $\tau$ and the inter-kick time $T$.
Red indicates that $\Lambda_{\mathrm{max}} > 0$; blue indicates that $\Lambda_{\mathrm{max}} < 0$.
The sign of $\Lambda_{\mathrm{max}}$ depends sensitively on $\tau$ and $T$ in the black region.  Here, $\sigma = 3$, $\lambda = 0.1$,  $A = 0.1$, and $\Phi = 1$.  For $t \in [-\tau, 0]$, we take $h(t) = (\theta(t), z(t)) = (0, t^2)$.
%\WOcomment{Bhargav, have I described the black region correctly?}
%\task{We need to add information about parameter values, our choice of $\Phi$, and our choice of history for these simulations.}
}
\label{fig:LinearDelay}
\end{figure}

\vspace{1.5ex}
\noindent
\textbf{Evidence for DIU.}
We argue that DIU emerges when $\tau > 0$.
Since no theory of rank one maps exists yet for infinite-dimensional delay systems, we compute the top Lyapunov exponent $\Lambda_{\mathrm{max}}$ numerically and use $\Lambda_{\mathrm{max}}$ as a DIU diagnostic.
That is, DIU is present if $\Lambda_{\mathrm{max}} > 0$ and absent if $\Lambda_{\mathrm{max}} < 0$.
We invoke the (\ref{li:excitability})--(\ref{li:external-forcing}) route to DIU and display our numerical evidence in Figure~\ref{fig:LinearDelay}.

(\ref{li:excitability}).
Delay induces excitability in system~\eqref{e:dlsf-kicked} because the stability of the limit cycle weakens as $\tau$ increases.
This is because the strength of stability of the zero solution to the scalar delay differential equation
\begin{equation*}
\dot{z} (t) = - \lambda z(t - \tau)
\end{equation*}
is determined by the (complex) solutions of the characteristic equation
\begin{equation*}
\gamma + \lambda e^{- \gamma \tau} = 0.
\end{equation*}
 
(\ref{li:intrinsic-shear}).
The angular velocity gradient parameter $\sigma$ in Eq.~\eqref{e:dlsf-intrinsic-thetaDot} quantifies shear in the unforced system.

(\ref{li:external-forcing}).
Here there exists a subtle difference between the delay ($\tau > 0$) case and the delay-free ($\tau = 0$) case.
In the delay-free case, the kick profile $\Phi$ needs to be nonconstant in order to create the $z$-variability that leads to stretching and folding of the phase space.
This requirement can be dropped in the delay case.
Here, the history $h$ that serves as initial data for system~\eqref{e:dlsf-kicked} can provide $z$-variability, so DIU can emerge even if $\Phi$ is a constant function.
Indeed, we demonstrate this in Figure~\ref{fig:LinearDelay}.

\vspace{1.5ex}
\noindent
\textbf{The delay-free case} $\bm{\tau = 0}$\textbf{.}
For the sake of completeness, we summarize what is known about this finite-dimensional system.
In the absence of delay, linear shear flow assumes the form
\begin{subequations}
\label{e:lsf-kicked}
\begin{align}
\dot{z} (t) &= - \lambda z (t) + A \Phi (\theta) \sum_{n=0}^{\infty} \delta (t - nT),
\label{e:lsf-kicked-z}
\\
\dot{\theta} &= 1 + \sigma z.
\label{e:lsf-kicked-theta}
\end{align}
\end{subequations}
This deceptively simple system has been studied by mathematical physicists such as Zaslavsky \cite{Zaslavsky1978simplest} and by dynamicists such as Lin and Young~\cite{Lin2008shear}.
Without forcing ($A = 0$), the set $\Omega = \{ (\theta, z) : z=0 \}$ is a limit cycle of system~\eqref{e:lsf-kicked} for all values of the contraction parameter $\lambda > 0$.
This limit cycle is weakly stable if $\lambda$ is small.
The dynamics of linear shear flow depend on the size and shape of the pulsatile kicks, as well as how trajectories relax to the limit cycle between kicks.

The angular velocity gradient parameter $\sigma$ quantifies the shear in system~\eqref{e:lsf-kicked}.
This parameter links with $A$ and $\lambda$ to form the {\itshape hyperbolicity factor}
\begin{equation*}
\frac{A \sigma}{\lambda} = \frac{(\text{kick amplitude}) \cdot (\text{shear})}{(\text{contraction})}.
\end{equation*}
The hyperbolicity factor governs the dynamics of system~\eqref{e:lsf-kicked}.
We express this through the behavior of the time-$T$ map $F_{T} : \mathbb{S}^{1} \times \mathbb{R} \to \mathbb{S}^{1} \times \mathbb{R}$ generated by linear shear flow (the result of one kick followed by one relaxation cycle).

If $\frac{A \sigma}{\lambda}$ is small, then the kicked system will quickly return to equilibrium after each kick (Figure~\ref{fig:sawtooth}a).
Mathematically, $F_{T}$ admits an attractor diffeomorphic to $\mbb{S}^{1}$.
However, stretch and fold geometry emerges when the hyperbolicity factor is large (Figure~\ref{fig:sawtooth}b).
Provided the kick profile $\Phi$ is not a constant function, each kick creates wave-like variability in the $z$-direction (Figure~\ref{fig:sawtooth}b, first image).
If the relaxation time $T$ is sizable, shear will then cause the waves to stretch and fold (Figure~\ref{fig:sawtooth}b, second image).
Chaotic behavior consequently emerges when the hyperbolicity factor is large.

Mathematically, $F_{T}$ admits a strange attractor for a set of $T$ values of positive Lebesgue measure, provided the hyperbolicity factor is sufficiently large and the kick profile $\Phi$ is not a constant function.
For such values of $T$, $F_{T}$ is genuinely nonuniformly hyperbolic.
In particular, it exhibits a positive Lyapunov exponent (the chaos is sustained in time).
The strange attractor supports a unique ergodic Sinai-Ruelle-Bowen measure.
This measure describes the asymptotic distribution of almost every orbit in the basin of attraction (the chaos is observable).
The Sinai-Ruelle-Bowen measure possesses strong statistical properties: The system obeys a dynamical version of the central limit theorem, correlations for H\"{o}lder observables decay exponentially, and a large deviations principle holds.

Found in~\cite{Wang2003strange}, precise statements and proofs of these rigorous results for linear shear flow rely on the theory of rank one maps.
Developed by Wang and Young~\cite{Wang2001strange, Wang2008toward, Wang2013dynamical}, rank one theory provides a platform for proving the existence of nonuniformly hyperbolic dynamics (sustained, observable chaos) within concrete systems of interest in the biochemical and physical sciences.

The proofs for linear shear flow involve verifying the hypotheses of the theory of rank one maps.
In particular, Wang and Young~\cite{Wang2003strange} analyze a certain `infinite relaxation' $T \to \infty$ limit of $F_{T}$.
This procedure produces the singular limit, a parametrized family of circle maps.
Rank one theory links the dynamics of the singular limit to those of $F_{T}$.
For linear shear flow, the hyperbolicity factor gives the amount of expansion in the singular limit.
Expansion in the singular limit links to sustained, observable chaos for $F_{T}$.

\section{Discussion}
\label{s:discussion}

\subsection{Summary}

In this paper, we have proposed a novel route by which delay can cause the onset of sustained temporal chaos in externally forced dynamical systems.
We call the resulting chaos \emph{delay-induced uncertainty}.
We have formulated a framework that may be used to diagnose the presence of DIU.
The DIU framework consists of three components:
(i) delay induces excitability in the unforced system, (ii) the unforced system possesses shear, and (iii) the external forcing interacts with the shear.
We have shown that this framework can transform limit cycles into strange attractors.

Guided by our DIU framework, we have demonstrated numerically that the Ultradian glucose-insulin model can produce DIU.
DIU emerges for both deterministic and random Dirac-$\delta$ forcing, as well as for meal-like carbohydrate input signals.
Our numerical experiments use the largest Lyapunov exponent $\Lambda_{\max}$ as a DIU diagnostic: $\Lambda_{\max} > 0$ indicates that DIU is present, while $\Lambda_{\max} < 0$ indicates that it is absent.

The Ultradian model is finite-dimensional because delay in this model assumes the form of a three-stage linear filter.
Because of this, the theory of rank one maps from smooth ergodic theory provides a rigorous dynamical profile of DIU for the Ultradian model.
No such theory exists for the infinite-dimensional flows generated by delay differential equations.
Nevertheless, we believe that the mechanisms in our DIU framework can cause infinite-dimensional delay systems to produce DIU.
We have supported this belief by showing that a delay linear shear flow model produces DIU.
This model allows us to illuminate both the geometry of DIU and how delay leads to DIU onset.

Because many physiological systems oscillate and possess meaningful sources of delay, we believe DIU is relevant throughout mathematical physiology.

\subsection{DIU in clinical medicine: Impact and outlook}

%Our results are directly applicable to ICU glucose care.
The existence of delay-induced uncertainty in the Ultradian model acknowledges and potentially explains the difficulties clinicians face when managing glucose in the ICU~\cite{Rabinstein2009hyperglycemia, Taylor1999prospective, Uyttendaele2020risk}.
%It may feel to clinicians that if they just tried harder or had a better protocol or even had a patient-specific protocol, they would choose better treatments.
Predicting the precise temporal evolution of a chaotic system is extremely challenging.  %, particularly for a clinician trying to monitor the chaotic trajectory of the glucose of a patient.
%If the process is chaotic, however, then predicting the precise evolution of the system becomes next to impossible.
Recently, artificial intelligence has been touted as a potential solution to problems in health care~\cite{Topol2019high}.
Our results imply that artificial intelligence may not be the solution because glucose-insulin metabolism is not just complicated but actually chaotic.
A better path might be to acknowledge the chaos and approach clinical challenges statistically, estimating expected variances so that decisions can be based on the distribution of likely outcomes~\cite{Albers2019simple}.
%Sampling more frequently can better track deviations and allow for interventions before dangerous glucose levels are reached.

Better matching of models to real ICU experience may help identify and avoid treatments that are likely to produce chaotic behavior.
For example, in a different context wherein the kicks represent impulsive insulin administration~\cite{Huang2012modeling, Song2014modeling} rather than nutrition, smaller, more frequent insulin doses may be better than larger, less frequent doses.
In the language of this paper, the former corresponds to lower kick amplitude and smaller inter-kick time.
Note that this hypothesis is consistent with analysis of impulsive insulin injection models~\cite{Huang2012modeling, Song2014modeling}.
In particular, Section~3.4 of~\cite{Huang2012modeling} describes how insulin injection dose and period should be selected so as to maintain ideal homeostatic glucose levels.
For this particular model, smaller doses of insulin administered more frequently are more efficient and effective.
Our work here may serve as a foundation for future work aimed at understanding uncertainty induced by delay.  

%Beyond the ICU, we believe DIU broadly impacts clinical decision-making.  
%When treating a patient, a clinician prescribes a medical intervention, or lack thereof, based on experience: the outcome of randomized clinical trials, observational studies, or heuristic understanding.
%In many cases clinicians are not afforded the luxury of quantitative bases for decisions.
%Patient treatment is based on results of trial-and-error: An intervention is prescribed, and, if it fails, a new intervention is prescribed based on the results.
%However, such clinical empiricism assumes reliability of the physiological process being treated. DIU is thus important for potentially elucidating inexplicable failures in treatments based on experience.
% 
Practical medicine introduces additional complexity we will analyze in future work, most notably \emph{nonstationarity}.
In this paper we have analyzed stationary models, meaning that model parameters do not change over time.
In the Ultradian context, we have therefore held the overall health state of the patient fixed, since the Ultradian model parameters represent overall health state~\cite{Albers2017personalized}.
This is often a reasonable assumption for glycemic management in the ICU since the timescales of medical interventions are relatively fast.
On longer timescales over which diseases such as diabetes affect patient health state, we anticipate that nonstationarity will be a common issue for DIU analysis in health-related physiologic settings.
%This is reasonable for glycemic management in the ICU, given the relatively fast treatment timescales.
%Over longer timescales, however, we anticipate that nonstationarity will be relevant for DIU analysis.
%We have set up a dichotomy in this paper:
%Assuming stationarity, does response to an external stimulus meaningfully depend on the state of the system (glucose level, plasma insulin level, interstitial insulin level in the Ultradian case) at the time the stimulus arrives?
%If no, we have reliability.
%If yes, we will observe DIU. 
%How does this dichotomy generalize to nonstationary systems?

The DIU dynamical profile has two facets.
Sustained temporal chaos renders the precise temporal evolution of \emph{individual trajectories} extremely difficult to predict.
However, dynamical systems may also be viewed through a statistical lens:
\emph{How do observables behave statistically?}
DIU induces sustained temporal chaos, yet observables exhibit strong statistical regularity.
Whether DIU is a benefit or a hazard therefore depends on context and point of view.
We briefly present two explicit examples to illustrate this point: glycemic management of enterally fed patients in an intensive care unit and glycemic self-management for a person with type-2 diabetes (T2DM).

Clinicians managing glucose levels in an ICU are not necessarily interested in the exact temporal evolution of individual glucose trajectories, but rather the boundaries of the achieved glucose levels.
Our simulations show that while the presence of DIU renders individual glucose trajectories unpredictable, it narrows the support of the invariant measure (projected onto the glucose variable).
This narrowing could in theory positively impact glycemic management.
Therefore, while clinicians seek to avoid chaos in general, it could be beneficial to design glycemic management protocols that deliberately induce DIU.

On the other hand, T2DM patients monitor their blood glucose levels at home in order to quantify activity impact.
Here, rational decision-making requires understanding the temporal evolution of individual glucose trajectories.
Glucose readings potentially become uninterpretable when DIU is present.
This would render self-management a nearly impossible task and suggests that DIU may negatively impact T2DM patients.

We emphasize that our thoughts on T2DM are speculative ideas for future work.
Mathematically, the T2DM scenario differs from ICU glycemic management in an important way.
While the latter involves glycemic oscillations, analyzing the former involves analyzing how the system returns to a \emph{fixed point} after nutritional kicks.
This suggests a dynamics problem:
Does DIU emerge when infinite-dimensional flows with fixed points are subjected to pulsatile forcing drives?
For finite-dimensional nonlinear flows, the theory of rank one maps has been used to analyze certain types of fixed points~\cite{Ott2008strange}.

Our results set the stage for important future work.
We have shown in this paper that the DIU phenomenon has the potential to confound medical decision-making and observational understanding.
% Said differently, we have identified the problem.
% Solving the problem begins with detection.
Next, we must link DIU to data-driven science.
How do we detect DIU using experimental data?
What mechanisms produce DIU, and do these mechanisms possess specific detectable statistical signatures?

%DIU phenomenology can be subtle and complex.
%Indeed, we have argued that DIU can have negative or positive impact, even in the same general physiological setting.
%When and how should we deliberately induce DIU?
%How can we mitigate negative impacts?

\subsection{Directions for future mathematical research}

Our work suggests several intriguing mathematical problems.
Can one develop a general theory of rank one maps that would apply to the flows on function spaces generated by delay differential equations?
There exist two approaches.
One could build a theory of rank one maps that applies directly to Banach spaces.
Alternatively, one could combine existing rank one theory with invariant manifold techniques, as Lu, Wang, and Young have done for supercritical Hopf bifurcations in certain parabolic PDEs~\cite{Lu2013strange}.
A general theory of rank one maps for delay differential equations would apply to DDE glucose-insulin models such as the two-delay model introduced in~\cite{Li2006modeling} and analyzed in~\cite{Li2007analysis}.
Proving that such systems exhibit DIU is an important next step given the broad applicability of DDE glucose-insulin modeling~\cite{Kissler2014determination, Wu2011physiological}.

The external forcing signals we have studied here have relatively simple structure.
What do more complex forcing signals produce in the DIU context?
For instance, what if the forcing assumes the form of a continuous-time stochastic process with a jump component, such as L\'{e}vy noise?

Nonstationary dynamical systems have received considerable attention over the past fifteen years.
Here, the dynamical model itself varies in time.
Glycemic management in the ICU fits naturally into this setting:
Ultradian model parameters likely drift over time as overall patient health state slowly changes.
Moreover, we do not know the statistics of this parameter variability.
(If we had this statistical information, we would be in the setting of random dynamical systems.)
Can one develop DIU theory for nonstationary dynamical systems?

It will be important to deepen the links between DIU and data-driven science.
Can one develop methods to deduce the presence of DIU using experimental data?
What impact does DIU have on data assimilation methods?
In particular, how does the fact that the Ultradian model exhibits DIU impact our ability to fit this model to ICU patient data?

Finally, do there exist additional routes to DIU?
What happens when delay acts on multiple timescales?

\section{Author contributions}

BK performed the numerical simulations.
All authors developed the analysis and wrote the paper.

\section{Funding}

This work has been partially supported by National Science Foundation grant DMS 1816315 (WO) and by National Institutes of Health grants R01 LM006910 (DA, GH), R01 LM012734 (DA, GH), and R01 GM117138 (BK, WO).

\section{Data availability}

Data and code will be available on {\ttfamily github}.
This information is available upon request as well.

\section{Nominal parameter values for the Ultradian model}

Table~\ref{table:model_parameters} lists the parameter values we have used for our DIU experiments.
See~\cite{Albers2017personalized} for information about this parameter set.

\begin{table}[ht]
%\begin{adjustwidth}{-2.25in}{0in} % Comment out/remove adjustwidth environment if table fits in text column.
\centering
\caption{Full list of parameters for the Ultradian glucose-insulin model \cite{Albers2017personalized}. Note that IIGU and IDGU denote insulin-independent glucose utilization and insulin-dependent glucose utilization, respectively.}
\begin{tabular}{|p{1.2cm}|p{2.7cm}|p{8cm}|}
\hline
\multicolumn{3}{|p{8cm}|}{\textbf{Ultradian model parameters}} \\ \hline
\Cline{2pt}{1-3} \hline
\Cline{2pt}{1-3}
Name & Nominal Value  & Meaning \\ \hline \hline
$V_p$  & $3$ L  & plasma volume  \\ \hline \hline
$V_i$  & $11$ L  & interstitial  volume \\ \hline \hline
$V_g$ & $10$ L  & glucose space \\ \hline \hline
$E$  & $0.2$ L min$^{-1}$ &   exchange rate for insulin between remote
and plasma compartments \\ \hline \hline
$t_p$  & $6$ min  & time constant for plasma insulin degradation (via
kidney and liver filtering) \\ \hline \hline
$t_i$  & $100$ min & time constant for remote insulin degradation (via
muscle and adipose tissue) \\ \hline \hline
$t_d$  & $12$ min  & delay between plasma insulin and glucose
production \\ \hline \hline
$R_m$  & $209$ mU min$^{-1}$  & linear constant affecting insulin secretion  \\ \hline \hline
$a_1$  & $6.6$ & exponential constant affecting insulin secretion \\ \hline \hline
$C_1$  & $300$ mg L$^{-1}$ & exponential constant affecting insulin secretion \\ \hline \hline
$C_2$  & $144$ mg L$^{-1}$  & exponential constant affecting IIGU \\ \hline \hline
$C_3$  & $100$ mg L$^{-1}$  & linear constant affecting IDGU \\ \hline \hline
$C_4$  & $80$ mU L$^{-1}$ & factor affecting IDGU \\ \hline \hline
$C_5$  & $26$ mU L$^{-1}$  & exponential constant affecting IDGU \\ \hline \hline
$U_b$  & $72$ mg min$^{-1}$  & linear constant affecting IIGU \\ \hline \hline
$U_0$  & $4$ mg min$^{-1}$ & linear constant affecting IDGU \\ \hline \hline
$U_m$  & $94$ mg min$^{-1}$  &  linear constant affecting IDGU \\ \hline \hline
$R_g$  & $180$ mg min$^{-1}$  & linear constant affecting IDGU \\ \hline \hline
$\alpha$  & $7.5$ & exponential constant affecting IDGU \\ \hline \hline
$\beta$  & $1.772$ & exponent affecting IDGU \\ \hline \hline
\end{tabular}
\label{table:model_parameters}
%\end{adjustwidth}
\end{table}

\FloatBarrier

\bibliographystyle{siamplain-chaos-rev2}
\bibliography{diu-v6-chaos-rev2}

\begin{thebibliography}{10}

\bibitem{Abarbanel1993analysis}
{\sc H.~D. Abarbanel, R.~Brown, J.~J. Sidorowich, and L.~S. Tsimring}, {\em The
  analysis of observed chaotic data in physical systems}, Reviews of modern
  physics, 65 (1993), p.~1331,
  \url{https://doi.org/10.1103/RevModPhys.65.1331}.

\bibitem{Albers2019simple}
{\sc D.~Albers, M.~Levine, M.~Sirlanci, and A.~Stuart}, {\em A simple modeling
  framework for prediction in the human glucose-insulin system}, arXiv preprint
  arXiv:1910.14193,  (2019).

\bibitem{Albers2012population}
{\sc D.~J. Albers, G.~Hripcsak, and M.~Schmidt}, {\em Population physiology:
  leveraging electronic health record data to understand human endocrine
  dynamics}, PLoS One, 7 (2012), p.~e48058,
  \url{https://doi.org/10.1371/journal.pone.0048058}.

\bibitem{Albers2017personalized}
{\sc D.~J. Albers, M.~Levine, B.~Gluckman, H.~Ginsberg, G.~Hripcsak, and
  L.~Mamykina}, {\em Personalized glucose forecasting for type 2 diabetes using
  data assimilation}, PLoS computational biology, 13 (2017),
  \url{https://doi.org/10.1371/journal.pcbi.1005232}.

\bibitem{Albers2018mechanistic}
{\sc D.~J. Albers, M.~E. Levine, A.~Stuart, L.~Mamykina, B.~Gluckman, and
  G.~Hripcsak}, {\em Mechanistic machine learning: how data assimilation
  leverages physiologic knowledge using {B}ayesian inference to forecast the
  future, infer the present, and phenotype}, Journal of the American Medical
  Informatics Association, 25 (2018), pp.~1392--1401,
  \url{https://doi.org/10.1093/jamia/ocy106}.

\bibitem{Bressloff2015frequency}
{\sc P.~C. Bressloff and B.~R. Karamched}, {\em A frequency-dependent decoding
  mechanism for axonal length sensing}, Frontiers in cellular neuroscience, 9
  (2015), p.~281, \url{https://doi.org/10.3389/fncel.2015.00281}.

\bibitem{Christini2002introduction}
{\sc D.~J. Christini and L.~Glass}, {\em Introduction: Mapping and control of
  complex cardiac arrhythmias}, Chaos: An Interdisciplinary Journal of
  Nonlinear Science, 12 (2002), pp.~732--739,
  \url{https://doi.org/10.1063/1.1504061}.

\bibitem{Claassen2013nonconvulsive}
{\sc J.~Claassen, A.~Perotte, D.~Albers, S.~Kleinberg, J.~M. Schmidt, B.~Tu,
  N.~Badjatia, H.~Lantigua, L.~J. Hirsch, S.~A. Mayer, et~al.}, {\em
  Nonconvulsive seizures after subarachnoid hemorrhage: multimodal detection
  and outcomes}, Annals of neurology, 74 (2013), pp.~53--64,
  \url{https://doi.org/10.1002/ana.23859}.

\bibitem{Drozdov1995model}
{\sc A.~Drozdov and H.~Khanina}, {\em A model for ultradian oscillations of
  insulin and glucose}, Mathematical and Computer Modelling, 22 (1995),
  pp.~23--38, \url{https://doi.org/10.1016/0895-7177(95)00108-E}.

\bibitem{Glass1988time}
{\sc L.~Glass, A.~Beuter, and D.~Larocque}, {\em Time delays, oscillations, and
  chaos in physiological control systems}, Mathematical Biosciences, 90 (1988),
  pp.~111--125, \url{https://doi.org/10.1016/0025-5564(88)90060-0}.

\bibitem{Glass1990chaos}
{\sc L.~Glass and C.~P. Malta}, {\em Chaos in multi-looped negative feedback
  systems.}, Journal of theoretical biology, 145 (1990), pp.~217--223,
  \url{https://doi.org/10.1016/s0022-5193(05)80127-4}.

\bibitem{Graham2020reduced}
{\sc E.~Graham, N.~Elhadad, and D.~Albers}, {\em Reduced model for female
  endocrine dynamics: Validation and functional variations}, arXiv preprint
  arXiv:2006.05034,  (2020).

\bibitem{Gupta2013transcriptional}
{\sc C.~Gupta, J.~López, W.~Ott, K.~Josić, and M.~Bennett}, {\em
  Transcriptional delay stabilizes bistable gene networks}, Physical Review
  Letters, 111 (2013), \url{https://doi.org/10.1103/PhysRevLett.111.058104}.

\bibitem{Hodgkin1952quantitative}
{\sc A.~L. Hodgkin and A.~F. Huxley}, {\em A quantitative description of
  membrane current and its application to conduction and excitation in nerve},
  The Journal of physiology, 117 (1952), p.~500,
  \url{https://doi.org/10.1113/jphysiol.1952.sp004764}.

\bibitem{Huang2012modeling}
{\sc M.~Huang, J.~Li, X.~Song, and H.~Guo}, {\em Modeling impulsive injections
  of insulin: {T}owards artificial pancreas}, SIAM Journal on Applied
  Mathematics, 72 (2012), pp.~1524--1548,
  \url{https://doi.org/10.1137/110860306}.

\bibitem{Josic2011stochastic}
{\sc K.~Josić, J.~López, W.~Ott, L.~Shiau, and M.~Bennett}, {\em Stochastic
  delay accelerates signaling in gene networks}, PLoS Computational Biology, 7
  (2011), \url{https://doi.org/10.1371/journal.pcbi.1002264}.

\bibitem{Karamched2015delayed}
{\sc B.~R. Karamched and P.~C. Bressloff}, {\em Delayed feedback model of
  axonal length sensing}, Biophysical journal, 108 (2015), pp.~2408--2419,
  \url{https://doi.org/10.1016/j.bpj.2015.03.055}.

\bibitem{Keener1998mathematical}
{\sc J.~Keener and J.~Sneyd}, {\em Mathematical physiology}, vol.~8 of
  Interdisciplinary Applied Mathematics, Springer-Verlag, New York, 1998.

\bibitem{Kissler2014determination}
{\sc S.~Kissler, C.~Cichowitz, S.~Sankaranarayanan, and D.~Bortz}, {\em
  Determination of personalized diabetes treatment plans using a two-delay
  model}, Journal of Theoretical Biology, 359 (2014), pp.~101--111,
  \url{https://doi.org/10.1016/j.jtbi.2014.06.005}.

\bibitem{Kyrychko2018enhancing}
{\sc Y.~N. Kyrychko and I.~B. Schwartz}, {\em Enhancing noise-induced switching
  times in systems with distributed delays}, Chaos, 28 (2018), pp.~063106, 9,
  \url{https://doi.org/10.1063/1.5034106}.

\bibitem{Li2007analysis}
{\sc J.~Li and Y.~Kuang}, {\em Analysis of a model of the glucose-insulin
  regulatory system with two delays}, SIAM Journal on Applied Mathematics, 67
  (2007), pp.~757--776, \url{https://doi.org/10.1137/050634001}.

\bibitem{Li2006modeling}
{\sc J.~Li, Y.~Kuang, and C.~Mason}, {\em Modeling the glucose-insulin
  regulatory system and ultradian insulin secretory oscillations with two
  explicit time delays}, Journal of Theoretical Biology, 242 (2006),
  pp.~722--735, \url{https://doi.org/10.1016/j.jtbi.2006.04.002}.

\bibitem{Li2004period}
{\sc T.-Y. Li and J.~A. Yorke}, {\em Period three implies chaos}, in The Theory
  of Chaotic Attractors, Springer, 2004, pp.~77--84.

\bibitem{Lin2008shear}
{\sc K.~K. Lin and L.-S. Young}, {\em Shear-induced chaos}, Nonlinearity, 21
  (2008), pp.~899--922, \url{https://doi.org/10.1088/0951-7715/21/5/002}.

\bibitem{Lu2013strange}
{\sc K.~Lu, Q.~Wang, and L.-S. Young}, {\em Strange attractors for periodically
  forced parabolic equations}, Mem. Amer. Math. Soc., 224 (2013), pp.~vi+85,
  \url{https://doi.org/10.1090/S0065-9266-2012-00669-1}.

\bibitem{Mackey1977oscillation}
{\sc M.~C. Mackey and L.~Glass}, {\em Oscillation and chaos in physiological
  control systems}, Science, 197 (1977), pp.~287--289,
  \url{https://doi.org/10.1126/science.267326}.

\bibitem{Ott2008strange}
{\sc W.~Ott}, {\em Strange attractors in periodically-kicked degenerate {H}opf
  bifurcations}, Communications in Mathematical Physics, 281 (2008),
  pp.~775--791, \url{https://doi.org/10.1007/s00220-008-0499-0}.

\bibitem{Ott2010from}
{\sc W.~Ott and M.~Stenlund}, {\em From limit cycles to strange attractors},
  Comm. Math. Phys., 296 (2010), pp.~215--249,
  \url{https://doi.org/10.1007/s00220-010-0994-y}.

\bibitem{Ott2005prevalence}
{\sc W.~Ott and J.~A. Yorke}, {\em Prevalence}, Bull. Amer. Math. Soc. (N.S.),
  42 (2005), pp.~263--290, \url{https://doi.org/10.1090/S0273-0979-05-01060-8}.

\bibitem{Rabinstein2009hyperglycemia}
{\sc A.~A. Rabinstein}, {\em Hyperglycemia in critical illness: lessons from
  nice-sugar}, Neurocritical care, 11 (2009), pp.~131--132.

\bibitem{Song2014modeling}
{\sc X.~Song, M.~Huang, and J.~Li}, {\em Modeling impulsive insulin delivery in
  insulin pump with time delays}, SIAM Journal on Applied Mathematics, 74
  (2014), pp.~1763--1785, \url{https://doi.org/10.1137/130933137}.

\bibitem{Song2005local}
{\sc Y.~Song and J.~Wei}, {\em Local hopf bifurcation and global periodic
  solutions in a delayed predator--prey system}, Journal of Mathematical
  Analysis and Applications, 301 (2005), pp.~1--21,
  \url{https://doi.org/10.1016/j.jmaa.2004.06.056}.

\bibitem{Sottile2018association}
{\sc P.~D. Sottile, D.~Albers, C.~Higgins, J.~Mckeehan, and M.~M. Moss}, {\em
  The association between ventilator dyssynchrony, delivered tidal volume, and
  sedation using a novel automated ventilator dyssynchrony detection
  algorithm}, Critical care medicine, 46 (2018), p.~e151,
  \url{https://doi.org/10.1097/CCM.0000000000002849}.

\bibitem{Stricker2008fast}
{\sc J.~Stricker, S.~Cookson, M.~Bennett, W.~Mather, L.~Tsimring, and
  J.~Hasty}, {\em A fast, robust and tunable synthetic gene oscillator},
  Nature, 456 (2008), pp.~516--519, \url{https://doi.org/10.1038/nature07389}.

\bibitem{Stroh2020estimating}
{\sc J.~Stroh, T.~Bennett, V.~Kheyfets, and D.~Albers}, {\em Estimating
  intracranial pressure via low-dimensional models: toward a practical tool for
  clinical decision support at multi-hour timescales}, bioRxiv,  (2020).

\bibitem{Sturis1991computer}
{\sc J.~Sturis, K.~Polonsky, E.~Mosekilde, and E.~Van~Cauter}, {\em Computer
  model for mechanisms underlying ultradian oscillations of insulin and
  glucose}, American Journal of Physiology - Endocrinology and Metabolism, 260
  (1991), pp.~E801--E809,
  \url{https://doi.org/10.1152/ajpendo.1991.260.5.E801}.

\bibitem{Taylor1999prospective}
{\sc S.~J. Taylor, S.~B. Fettes, C.~Jewkes, and R.~J. Nelson}, {\em
  Prospective, randomized, controlled trial to determine the effect of early
  enhanced enteral nutrition on clinical outcome in mechanically ventilated
  patients suffering head injury}, Critical care medicine, 27 (1999),
  pp.~2525--2531.

\bibitem{Topol2019high}
{\sc E.~J. Topol}, {\em High-performance medicine: the convergence of human and
  artificial intelligence}, Nature medicine, 25 (2019), pp.~44--56.

\bibitem{Urteaga2019multi}
{\sc I.~Urteaga, T.~Bertin, T.~M. Hardy, D.~J. Albers, and N.~Elhadad}, {\em
  Multi-task gaussian processes and dilated convolutional networks for
  reconstruction of reproductive hormonal dynamics}, arXiv preprint
  arXiv:1908.10226,  (2019).

\bibitem{Uyttendaele2020risk}
{\sc V.~Uyttendaele, J.~Knopp, G.~Shaw, T.~Desaive, and J.~Chase}, {\em Risk
  and reward: {E}xtending stochastic glycaemic control intervals to reduce
  workload}, BioMedical Engineering Online, 19 (2020),
  \url{https://doi.org/10.1186/s12938-020-00771-6}.

\bibitem{Wang2001strange}
{\sc Q.~Wang and L.-S. Young}, {\em Strange attractors with one direction of
  instability}, Comm. Math. Phys., 218 (2001), pp.~1--97,
  \url{https://doi.org/10.1007/s002200100379}.

\bibitem{Wang2003strange}
{\sc Q.~Wang and L.-S. Young}, {\em Strange attractors in periodically-kicked
  limit cycles and {H}opf bifurcations}, Comm. Math. Phys., 240 (2003),
  pp.~509--529, \url{https://doi.org/10.1007/s00220-003-0902-9}.

\bibitem{Wang2008toward}
{\sc Q.~Wang and L.-S. Young}, {\em Toward a theory of rank one attractors},
  Ann. of Math. (2), 167 (2008), pp.~349--480,
  \url{https://doi.org/10.4007/annals.2008.167.349}.

\bibitem{Wang2013dynamical}
{\sc Q.~Wang and L.-S. Young}, {\em Dynamical profile of a class of rank-one
  attractors}, Ergodic Theory Dynam. Systems, 33 (2013), pp.~1221--1264,
  \url{https://doi.org/10.1017/S014338571200020X}.

\bibitem{Wei2005hopf}
{\sc J.~Wei and M.~Y. Li}, {\em Hopf bifurcation analysis in a delayed
  nicholson blowflies equation}, Nonlinear Analysis: Theory, Methods \&
  Applications, 60 (2005), pp.~1351--1367,
  \url{https://doi.org/10.1016/j.na.2003.04.002}.

\bibitem{Wu2011physiological}
{\sc Z.~Wu, C.-K. Chui, G.-S. Hong, and S.~Chang}, {\em Physiological analysis
  on oscillatory behavior of glucose-insulin regulation by model with delays},
  Journal of Theoretical Biology, 280 (2011), pp.~1--9,
  \url{https://doi.org/10.1016/j.jtbi.2011.03.032}.

\bibitem{Xu2017theory}
{\sc B.~Xu and P.~C. Bressloff}, {\em A theory of synchrony for active
  compartments with delays coupled through bulk diffusion}, Physica D:
  Nonlinear Phenomena, 341 (2017), pp.~45--59,
  \url{https://doi.org/10.1016/j.physd.2016.10.001}.

\bibitem{Yan2006hopf}
{\sc X.-P. Yan and W.-T. Li}, {\em Hopf bifurcation and global periodic
  solutions in a delayed predator--prey system}, Applied Mathematics and
  Computation, 177 (2006), pp.~427--445,
  \url{https://doi.org/10.1016/j.amc.2005.11.020}.

\bibitem{Zaslavsky1978simplest}
{\sc G.~M. Zaslavsky}, {\em The simplest case of a strange attractor}, Phys.
  Lett. A, 69 (1978/79), pp.~145--147,
  \url{https://doi.org/10.1016/0375-9601(78)90195-0}.

\bibitem{Zenker2007inverse}
{\sc S.~Zenker, J.~Rubin, and G.~Clermont}, {\em From inverse problems in
  mathematical physiology to quantitative differential diagnoses}, PLoS Comput
  Biol, 3 (2007), p.~e204, \url{https://doi.org/10.1371/journal.pcbi.1007155}.

\end{thebibliography}

\end{document}